\documentclass[11pt]{article}

	\newcommand{\mainTitle}{Extended double shuffle relations and the generating function of triple zeta values of any fixed weight}

	\newcommand{\authorName}{Tomoya Machide}
	\newcommand{\organizationName}{Kinki University}
	\newcommand{\departmentName}{Interdisciplinary Graduate School of Science and Engineering}
	
	\newcommand{\placeAddress}{3-4-1 Kowakae, Higashi-Osaka, Osaka 577-8502, Japan}
	\newcommand{\emailAddress}{E-mail: machide.t@gmail.com}
	\newcommand{\sectOne}{Introduction and statement of results}
	\newcommand{\sectTwo}{Generating function of multiple polylogarithm}
	\newcommand{\sSectTwoO}{Formulas corresponding to shuffle and harmonic relations}
	\newcommand{\sSectTwoT}{Asymptotic properties}
	\newcommand{\sectThree}{Formulas via extended double shuffle relations}%Formulas derived from extended double shuffle relations
	\newcommand{\sectFour}{Parameterized and weighted sum formulas}%Parameterized sum formula
	\newcommand{\sectFive}{Restricted sum formulas}%Restricted sum formulas
\usepackage{geometry}               
\usepackage{ascmac}
\usepackage{amsmath}
\usepackage{amssymb}
\usepackage{amsthm}
\usepackage[hypertex, colorlinks=true, bookmarks=true]{hyperref}

	\DeclareMathOperator*{\OPlus}{\bigoplus}
	\newcommand{\nbk}[3]{#1#3#2}		
	\newcommand{\bgbk}[3]{\bigl{#1}#3\bigr{#2}}	
	\newcommand{\Bgbk}[3]{\Bigl{#1}#3\Bigr{#2}}			
	\newcommand{\bggbk}[3]{\biggl{#1}#3\biggr{#2}}			
	\newcommand{\Bggbk}[3]{\Biggl{#1}#3\Biggr{#2}}
	\newcommand{\autobk}[3]{\left#1#3\right#2}
	\newcommand{\mcbk}[4][?]{\ifx n#1\nbk{#2}{#3}{#4}\else\ifx b#1\bgbk{#2}{#3}{#4}\else\ifx B#1\Bgbk{#2}{#3}{#4}\else\ifx g#1\bggbk{#2}{#3}{#4}\else\ifx G#1\Bggbk{#2}{#3}{#4}\else\ifx a#1\autobk{#2}{#3}{#4}\else#4\fi\fi\fi\fi\fi\fi}
	\newcommand{\nsgsb}[1]{#1}		
	\newcommand{\bgsgsb}[1]{\big{#1}}	
	\newcommand{\Bgsgsb}[1]{\Big{#1}}			
	\newcommand{\bggsgsb}[1]{\bigg{#1}}			
	\newcommand{\Bggsgsb}[1]{\Bigg{#1}}
	\newcommand{\mcsgsb}[2][?]{\ifx n#1\nsgsb{#2}\else\ifx b#1\bgsgsb{#2}\else\ifx B#1\Bgsgsb{#2}\else\ifx g#1\bggsgsb{#2}\else\ifx G#1\Bggsgsb{#2}\else#2\fi\fi\fi\fi\fi}
	\newcommand{\myEqSpace}{\,}
	\newlength{\myEqSpaceLen}
	\newcommand{\setZ}{\mathbb{Z}} 
	\newcommand{\setQ}{\mathbb{Q}}	
	\newcommand{\setR}{\mathbb{R}}

	\newcommand{\gpSym}[2][?]{\ifx g#1\mathfrak{S}_{#2} \else S_{#2} \fi}
	\newcommand{\gpAlt}[2][?]{\ifx g#1 \mathfrak{A}_{#2} \else A_{#2} \fi }
	\newcommand{\gpKleinF}[1][?]{\ifx g#1 \mathfrak{V} \else V \fi }
	\newcommand{\gpCyc}[2][?]{\ifx g#1 \mathfrak{C}_{#2} \else C_{#2}\fi}	
	
	\newcommand{\lnA}[1][]{&  &}
	\newcommand{\lnP}[1]{\myEqSpace#1\myEqSpace}
	\newcommand{\lnAP}[2][]{& #2 &}
	\newcommand{\lnAH}[1][\nonumber]{#1\\ & &}
	\newcommand{\lnAHP}[2][\nonumber]{#1 \\ & #2 &}
	\newcommand{\lnAHs}[2][\nonumber]{#1 \\ & &\hspace{#2pt}}

	\newcommand{\refH}[2]{#1\ref{#2}}
	\newcommand{\refF}[3][\ ]{#2#1#3}
	\newcommand{\refEq}[1]{(\ref{#1})}	
		
	\newcommand{\refSect}[2][?]{\refH{\S}{#2}}
	
	\newcommand{\bkR}[2][a]{\mcbk[#1]{(}{)}{#2}}
	\newcommand{\bkS}[2][a]{\mcbk[#1]{[}{]}{#2}}
	\newcommand{\bkB}[2][a]{\mcbk[#1]{\{}{\}}{#2}}
	\newcommand{\bkA}[2][a]{\mcbk[#1]{\langle}{\rangle}{#2}}

	\newcommand{\tmR}[4][?]{\bkR[#1]{#2} \ifx 0#4 \else \ifx 1#4 #3 \else {#3}^{#4} \fi \fi}	
	\newcommand{\tmS}[4][?]{\bkS[#1]{#2} \ifx 0#4 \else \ifx 1#4 #3 \else {#3}^{#4} \fi \fi}
	\newcommand{\tmM}[4][?]{\bkB[#1]{#2} \ifx 0#4 \else \ifx 1#4 #3 \else {#3}^{#4} \fi \fi}

	\newcommand{\alp}{\alpha}

	\newcommand{\ep}{\varepsilon}
	\newcommand{\sig}{\sigma}

	\newcommand{\mo}{(-1)}

	\newcommand{\mVert}[1][n]{{\ \mcsgsb[#1]{\vert}\ }}
	
	\newcommand{\opF}[3][?]{\ifx s#1 #2/#3 \else \ifx d#1 \dfrac{#2}{#3} \else \frac{#2}{#3} \fi\fi}

	\newcommand{\pw}[3][?]{\ifx 0#3 1  \else\ifx 1#3 \bkR[#1]{#2}\else \ifx?#1\bkR[#1]{#2}^{#3}\else {\bkR[#1]{#2}}^{#3} \fi\fi\fi}
	\newcommand{\id}[3][?]{#2_{#3}}
	\newcommand{\ip}[4][?]{\ifx ?#1 \ifx 1#4 {#2}_{#3} \else {#2}_{#3}^{#4} \fi \else \ifx n#1\ifx 1#3 {(#2_{#3})} \else (#2_{#3})^{#4} \fi \else #1\fi\fi}
	\newcommand{\pwR}[3][a]{ \pw[#1]{#2}{#3} }
	\newcommand{\pwB}[3][a]{\ifx 0#3 1  \else\ifx 1#3 \bkB[#1]{#2}\else \ifx?#1\bkB[#1]{#2}^{#3}\else {\bkB[#1]{#2}}^{#3} \fi\fi\fi}
					
	\newcommand{\sqAc}[3][]{#2#1\cdots#1#3}
		
		\newcommand{\sqSm}[1][+]{\sqAc[#1]}
		
	\newcommand{\nFc}[3][n]{#2\bkR[#1]{#3}}
	
	\newcommand{\idFc}[4][n]{\id{#2}{#3}\bkR[#1]{#4}}
	\newcommand{\pwFc}[4][n]{\pw{#2}{#3}\bkR[#1]{#4}}
	\newcommand{\ipFc}[5][n]{\ip{#2}{#3}{#4}\bkR[#1]{#5}}
		\newcommand{\Fc}{\nFc}
	
	\newcommand{\tpT}[3][a]{ {#2}\atop \bkR[#1]{#3} }
	\newcommand{\tpTh}[4][a]{{#2}\atop \bkR[#1]{{#3\atop#4}} }

	\newcommand{\nSm}[2][?]{\ifx l#1 \sum\limits_{#2} \else \sum_{#2} \fi}
	\newcommand{\nSmT}[3][?]{\ifx l#1 \sum\limits_{#2}^{#3} \else \sum_{#2}^{#3} \fi}	
	\newcommand{\pSm}[2][?]{\ifx t#1 \sum_{#2}^{\prime} \else \sideset{}{^\prime}\sum_{#2} \fi}
	\newcommand{\pSmT}[3][?]{\ifx t#1 \sum_{#2}^{\prime#3} \else \sideset{}{^\prime}\sum_{#2}^{#3} \fi}	
	\newcommand{\tpTSm}[3][?]{\nSm[#1]{\tpT{#2}{#3}}}
	
	\newcommand{\tpThSm}[4][?]{\nSm[#1]{\tpTh{#2}{#3}{#4}}}
	
		\newcommand{\Sm}{\nSm}
		\newcommand{\SmT}{\nSmT}
		\newcommand{\tpSm}{\tpTSm}
		
	\newcommand{\nPd}[2][?]{\ifx l#1 \prod\limits_{#2} \else \prod_{#2} \fi}
	\newcommand{\nPdT}[3][?]{\ifx l#1 \prod\limits_{#2}^{#3} \else \prod_{#2}^{#3} \fi}

		\newcommand{\PdT}{\nPdT}

	\newcommand{\nOPs}[2][?]{\ifx l#1 \OPlus\limits_{#2} \else \OPlus_{#2} \fi}
	\newcommand{\nOPsT}[3][?]{\ifx l#1 \OPlus\limits_{#2}^{#3} \else \OPlus_{#2}^{#3} \fi}	
	\newcommand{\pOPs}[2][?]{\ifx t#1 \OPlus_{#2}^{\prime} \else \sideset{}{^\prime}\OPlus_{#2} \fi}
	\newcommand{\pOPsT}[3][?]{\ifx t#1 \OPlus_{#2}^{\prime#3} \else \sideset{}{^\prime}\OPlus_{#2}^{#3} \fi}

	\newcommand{\nIs}[2][?]{\ifx l#1 \bigcap\limits_{#2} \else \bigcap_{#2} \fi}
	\newcommand{\nIsT}[3][?]{\ifx l#1 \bigcap\limits_{#2}^{#3} \else \bigcap_{#2}^{#3} \fi}	
	\newcommand{\pIs}[2][?]{\ifx t#1 \bigcap_{#2}^{\prime} \else \sideset{}{^\prime}\bigcap_{#2} \fi}
	\newcommand{\pIsT}[3][?]{\ifx t#1 \bigcap_{#2}^{\prime#3} \else \sideset{}{^\prime}\bigcap_{#2}^{#3} \fi}

	\newcommand{\nUn}[2][?]{\ifx l#1 \bigcup\limits_{#2} \else \bigcup_{#2} \fi}
	\newcommand{\nUnT}[3][?]{\ifx l#1 \bigcup\limits_{#2}^{#3} \else \bigcup_{#2}^{#3} \fi}	
	\newcommand{\pUn}[2][?]{\ifx t#1 \bigcup_{#2}^{\prime} \else \sideset{}{^\prime}\bigcup_{#2} \fi}
	\newcommand{\pUnT}[3][?]{\ifx t#1 \bigcup_{#2}^{\prime#3} \else \sideset{}{^\prime}\bigcup_{#2}^{#3} \fi}

	\newcommand{\nLim}[2][?]{\ifx l#1 \lim\limits_{#2} \else \lim_{#2} \fi}
		
	\newcommand{\glcondEnvLineHead}[1]{ \ifx*#1 \begin{eqnarray*} \else \begin{eqnarray}  \label{#1} \fi }
	\newcommand{\glcondEnvLineTail}[1]{ \ifx*#1 \end{eqnarray*} \else \end{eqnarray} \fi }
	\newcommand{\glcondDis}[1]{\ifx d#1 \displaystyle \fi}

		\newcommand{\envHLineT}[3][*]{ \glcondEnvLineHead{#1} #2&=&#3\glcondEnvLineTail{#1} }
		\newcommand{\envHLineTDef}[3][*]{ \glcondEnvLineHead{#1} #2&:=&#3\glcondEnvLineTail{#1} }
			\newcommand{\envHLine}{\envHLineT}
			\newcommand{\envHLineDef}{\envHLineTDef}
		\newcommand{\envHLineCE}[9][*]{\glcondEnvLineHead{#1} #2&=&#3\\#4&=&#5\nonumber\\#6&=&#7\nonumber\\#8&=&#9\nonumber\glcondEnvLineTail{#1}}
		
		\newcommand{\envPLine}[2][*]{\glcondEnvLineHead{#1} #2\glcondEnvLineTail{#1}}
		\newcommand{\envMO}[2][*]{$\ifx d#1 \displaystyle \fi#2$}
		\newcommand{\envMT}[3][*]{$\ifx d#1 \displaystyle \fi#2=#3$}
		\newcommand{\envMTDef}[3][*]{$\ifx d#1 \displaystyle \fi#2:=#3$}
			\newcommand{\envM}{\envMT}
			\newcommand{\envMDef}{\envMTDef}
		\newcommand{\envMTh}[4][*]{$\ifx d#1 \displaystyle \fi#2=#3=#4$}
		\newcommand{\envMF}[5][*]{$\ifx d#1 \displaystyle \fi#2=#3=#4=#5$}
		
	\newcommand{\envMLineT}[3][*]{ \ifx*#1 \begin{multline*} #2\lnP{=}#3\end{multline*} \else \begin{multline} \label{#1} #2\lnP{=}#3\end{multline} \fi }

	\newcommand{\lcparaCase}{\vspace{3pt}}
	\newcommand{\envCaseT}[3][?]{\begin{cases} \glcondDis{#1}#2,\lcparaCase\\\glcondDis{#1}#3\end{cases}}

	\newcommand{\envSMatT}[3][a]{ \autobk{(}{)}{\begin{smallmatrix}#2\\#3\end{smallmatrix}} }

	\newcommand{\abs}[1]{\left | #1 \right |  }		
	\newcommand{\sbLandau}[2][a]{\Fc[#1]{O}{#2}}

	\newcommand{\lgg}[2][?]{\Fc[#1]{\log}{#2}}

	\newcommand{\fcZeta}[1]{\zeta(#1)}
		\newcommand{\fcZ}{\fcZeta}		
	\newcommand{\lgP}[3][n]{\idFc[#1]{Li}{#2}{#3}}
	
	\newcommand{\envMLineTPt}[4][*]{ \ifx*#1 \begin{multline*} #3\lnP{#2}#4\end{multline*} \else \begin{multline} \label{#1} #3\lnP{#2}#4\end{multline} \fi }

	\newcommand{\glcondEnvLineTailPd}[1]{.\ifx*#1 \end{eqnarray*} \else \end{eqnarray} \fi }
	\newcommand{\glcondEnvLineTailCm}[1]{,\ifx*#1 \end{eqnarray*} \else \end{eqnarray} \fi }
	\newcommand{\envProof}[2][?]{ \par\mbox{}\vspace{-5pt}\\ \ifx?#1\emph{Proof.}\else\emph{Proof of #1.}\fi \ #2 \hfill $\Box$\\ \par}
	
		\newcommand{\envLineTPd}[3][*]{ \glcondEnvLineHead{#1} & &#2\\&=&#3\nonumber \glcondEnvLineTailPd{#1} }
		
		\newcommand{\envLineTCm}[3][*]{ \glcondEnvLineHead{#1} & &#2\\&=&#3\nonumber \glcondEnvLineTailCm{#1} }
		
			\newcommand{\envLinePd}{\envLineTPd}
			
			\newcommand{\envLineCm}{\envLineTCm}
			
		\newcommand{\envLineFPd}[5][*]{ \glcondEnvLineHead{#1} & &#2\\&=&#3\nonumber \\&=&#4\nonumber \\&=&#5\nonumber \glcondEnvLineTailPd{#1} }
		\newcommand{\envLineFCm}[5][*]{ \glcondEnvLineHead{#1} & &#2\\&=&#3\nonumber \\&=&#4\nonumber \\&=&#5\nonumber \glcondEnvLineTailCm{#1} }
		\newcommand{\envHLineTPd}[3][*]{ \glcondEnvLineHead{#1} #2&=&#3\glcondEnvLineTailPd{#1} }
		\newcommand{\envHLineTDefPd}[3][*]{ \glcondEnvLineHead{#1} #2&:=&#3\glcondEnvLineTailPd{#1} }
		\newcommand{\envHLineTCm}[3][*]{ \glcondEnvLineHead{#1} #2&=&#3\glcondEnvLineTailCm{#1} }
		
			\newcommand{\envHLinePd}{\envHLineTPd}
			\newcommand{\envHLineDefPd}{\envHLineTDefPd}
			\newcommand{\envHLineCm}{\envHLineTCm}
			
		\newcommand{\envHLineThPd}[4][*]{ \glcondEnvLineHead{#1} #2&=&#3\\&=&#4\nonumber\glcondEnvLineTailPd{#1}	 }

		\newcommand{\envHLineFPd}[5][*]{ \glcondEnvLineHead{#1} #2&=&#3\\&=&#4\nonumber \\&=&#5\nonumber \glcondEnvLineTailPd{#1} }
		
		\newcommand{\envHLineFiCm}[6][*]{ \glcondEnvLineHead{#1} #2&=&#3\\&=&#4\nonumber \\&=&#5\nonumber  \\&=&#6\nonumber\glcondEnvLineTailCm{#1}}
		\newcommand{\envHLineCFPd}[5][*]{\glcondEnvLineHead{#1} #2&=&#3,\\#4&=&#5\nonumber\glcondEnvLineTailPd{#1}}
		\newcommand{\envHLineCFDefPd}[5][*]{\glcondEnvLineHead{#1} #2&:=&#3,\\#4&:=&#5\nonumber\glcondEnvLineTailPd{#1}}
		
		\newcommand{\envHLineCFCmDef}[5][*]{\glcondEnvLineHead{#1} #2&:=&#3,\\#4&:=&#5\nonumber\glcondEnvLineTailCm{#1}}
		\newcommand{\envHLineCFNmePd}[5][?]{\begin{eqnarray} #2&=&#3,\\#4&=&#5 \glcondEnvLineTailPd{?} }
		\newcommand{\envHLineCFDefNmePd}[5][?]{\begin{eqnarray} #2&:=&#3,\\#4&:=&#5 \glcondEnvLineTailPd{?} }
		\newcommand{\envHLineCFCmNme}[5][?]{\begin{eqnarray} #2&=&#3,\\#4&=&#5 \glcondEnvLineTailCm{?} }
		\newcommand{\envHLineCFCmDefNme}[5][?]{\begin{eqnarray} #2&:=&#3,\\#4&:=&#5 \glcondEnvLineTailCm{?} }
		\newcommand{\envHLineCSPd}[7][*]{\glcondEnvLineHead{#1} #2&=&#3,\\#4&=&#5,\nonumber\\#6&=&#7\nonumber\glcondEnvLineTailPd{#1}}
		\newcommand{\envHLineCSDefPd}[7][*]{\glcondEnvLineHead{#1} #2&:=&#3,\\#4&:=&#5,\nonumber\\#6&:=&#7\nonumber\glcondEnvLineTailPd{#1}}
		\newcommand{\envHLineCSCm}[7][*]{\glcondEnvLineHead{#1} #2&=&#3,\\#4&=&#5,\nonumber\\#6&=&#7\nonumber\glcondEnvLineTailCm{#1}}
			
		\newcommand{\envHLineCSNmePd}[7][*]{\begin{eqnarray} #2&=&#3,\\#4&=&#5,\\#6&=&#7\glcondEnvLineTailPd{?}}
		\newcommand{\envHLineCSDefNmePd}[7][*]{\begin{eqnarray} #2&:=&#3,\\#4&:=&#5,\\#6&:=&#7\glcondEnvLineTailPd{?}}
		\newcommand{\envHLineCSCmNme}[7][*]{\begin{eqnarray} #2&=&#3,\\#4&=&#5,\\#6&=&#7\glcondEnvLineTailCm{?}}
		\newcommand{\envHLineCSCmDefNme}[7][*]{\begin{eqnarray} #2&:=&#3,\\#4&:=&#5,\\#6&:=&#7\glcondEnvLineTailCm{?}}
		\newcommand{\envHLineCEPd}[9][*]{\glcondEnvLineHead{#1} #2&=&#3,\\#4&=&#5,\nonumber\\#6&=&#7,\nonumber\\#8&=&#9\nonumber\glcondEnvLineTailPd{#1}}

		\newcommand{\envHLineCENmePd}[9][*]{\begin{eqnarray} #2&=&#3,\\#4&=&#5,\\#6&=&#7,\\#8&=&#9\glcondEnvLineTailPd{?}}
		\newcommand{\envHLineCEDefNmePd}[9][*]{\begin{eqnarray} #2&:=&#3,\\#4&:=&#5,\\#6&:=&#7,\\#8&:=&#9\glcondEnvLineTailPd{?}}
		\newcommand{\envHLineCECmNme}[9][*]{\begin{eqnarray} #2&=&#3,\\#4&=&#5,\\#6&=&#7,\\#8&=&#9\glcondEnvLineTailCm{?}}
		\newcommand{\envHLineCECmDefNme}[9][*]{\begin{eqnarray} #2&:=&#3,\\#4&:=&#5,\\#6&:=&#7,\\#8&:=&#9\glcondEnvLineTailCm{?}}
		\newcommand{\envPLinePd}[2][*]{\glcondEnvLineHead{#1} #2\glcondEnvLineTailPd{#1}}
		\newcommand{\envPLineCm}[2][*]{\glcondEnvLineHead{#1} #2\glcondEnvLineTailCm{#1}}
		\newcommand{\envOTLinePd}[4][*]{\glcondEnvLineHead{#1} #2\lnAP{=}#3\lnP{=}#4.\glcondEnvLineTail{#1}}
		\newcommand{\envOTLineCm}[4][*]{\glcondEnvLineHead{#1} #2\lnAP{=}#3\lnP{=}#4,\glcondEnvLineTail{#1}}

		\newcommand{\envMOPd}[2][*]{$\ifx d#1 \displaystyle \fi#2.$}
		\newcommand{\envMOCm}[2][*]{$\ifx d#1 \displaystyle \fi#2,$}
		\newcommand{\envMTPd}[3][*]{$\ifx d#1 \displaystyle \fi#2=#3.$}
		\newcommand{\envMTCm}[3][*]{$\ifx d#1 \displaystyle \fi#2=#3,$}
		\newcommand{\envMTDefPd}[3][*]{$\ifx d#1 \displaystyle \fi#2:=#3.$}
		\newcommand{\envMTDefCm}[3][*]{$\ifx d#1 \displaystyle \fi#2:=#3,$}
			\newcommand{\envMPd}{\envMTPd}
			\newcommand{\envMCm}{\envMTCm}

		\newcommand{\envHLineCFCmNm}[5][*]{ \begin{equation}\begin{split} \ifx*#1 \text{[ERROR;need label name]} \else \label{#1} \fi #2&\lnP{=}#3,\\#4&\lnP{=}#5, \end{split}\end{equation} }
		\newcommand{\envHLineCFNm}[5][*]{ \begin{equation}\begin{split} \ifx*#1 \text{[ERROR;need label name]} \else \label{#1} \fi #2&\lnP{=}#3\\#4&\lnP{=}#5, \end{split}\end{equation} }
		\newcommand{\envHLineCFNmPd}[5][*]{ \begin{equation}\begin{split} \ifx*#1 \text{[ERROR;need label name]} \else \label{#1} \fi #2&\lnP{=}#3,\\#4&\lnP{=}#5. \end{split}\end{equation} }
		\newcommand{\envHLineCFCmDefNm}[5][*]{ \begin{equation}\begin{split} \ifx*#1 \text{[ERROR;need label name]} \else \label{#1} \fi #2&\lnP{:=}#3,\\#4&\lnP{:=}#5, \end{split}\end{equation} }
		\newcommand{\envHLineCFDefNm}[5][*]{ \begin{equation}\begin{split} \ifx*#1 \text{[ERROR;need label name]} \else \label{#1} \fi #2&\lnP{:=}#3\\#4&\lnP{:=}#5, \end{split}\end{equation} }
		\newcommand{\envHLineCFDefNmPd}[5][*]{ \begin{equation}\begin{split} \ifx*#1 \text{[ERROR;need label name]} \else \label{#1} \fi #2&\lnP{:=}#3,\\#4&\lnP{:=}#5. \end{split}\end{equation} }
		\newcommand{\envHLineCSCmNm}[7][*]{ \begin{equation}\begin{split} \ifx*#1 \text{[ERROR;need label name]} \else \label{#1} \fi #2&\lnP{=}#3,\\#4&\lnP{=}#5,\\#6&\lnP{=}#7 \end{split}\end{equation} }
		\newcommand{\envHLineCSNm}[7][*]{ \begin{equation}\begin{split} \ifx*#1 \text{[ERROR;need label name]} \else \label{#1} \fi #2&\lnP{=}#3\\#4&\lnP{=}#5\\#6&\lnP{=}#7 \end{split}\end{equation} }
		\newcommand{\envHLineCSNmPd}[7][*]{ \begin{equation}\begin{split} \ifx*#1 \text{[ERROR;need label name]} \else \label{#1} \fi #2&\lnP{=}#3,\\#4&\lnP{=}#5,\\#6&\lnP{=}#7. \end{split}\end{equation} }
		\newcommand{\envHLineCSCmDefNm}[7][*]{ \begin{equation}\begin{split} \ifx*#1 \text{[ERROR;need label name]} \else \label{#1} \fi #2&\lnP{:=}#3,\\#4&\lnP{:=}#5,\\#6&\lnP{:=}#7 \end{split}\end{equation} }
		\newcommand{\envHLineCSDefNm}[7][*]{ \begin{equation}\begin{split} \ifx*#1 \text{[ERROR;need label name]} \else \label{#1} \fi #2&\lnP{:=}#3\\#4&\lnP{:=}#5\\#6&\lnP{:=}#7 \end{split}\end{equation} }
		\newcommand{\envHLineCSDefNmPd}[7][*]{ \begin{equation}\begin{split} \ifx*#1 \text{[ERROR;need label name]} \else \label{#1} \fi #2&\lnP{:=}#3,\\#4&\lnP{:=}#5,\\#6&\lnP{:=}#7. \end{split}\end{equation} }
		\newcommand{\envHLineCECmNm}[9][*]{ \begin{equation}\begin{split} \ifx*#1 \text{[ERROR;need label name]} \else \label{#1} \fi #2&\lnP{=}#3,\\#4&\lnP{=}#5,\\#6&\lnP{=}#7,\\#8&\lnP{=}#9,  \end{split}\end{equation} }
		\newcommand{\envHLineCENm}[9][*]{ \begin{equation}\begin{split} \ifx*#1 \text{[ERROR;need label name]} \else \label{#1} \fi #2&\lnP{=}#3\\#4&\lnP{=}#5\\#6&\lnP{=}#7\\#8&\lnP{=}#9  \end{split}\end{equation} }
		\newcommand{\envHLineCENmPd}[9][*]{ \begin{equation}\begin{split} \ifx*#1 \text{[ERROR;need label name]} \else \label{#1} \fi #2&\lnP{=}#3,\\#4&\lnP{=}#5,\\#6&\lnP{=}#7,\\#8&\lnP{=}#9.  \end{split}\end{equation} }
		\newcommand{\envHLineCECmDefNm}[9][*]{ \begin{equation}\begin{split} \ifx*#1 \text{[ERROR;need label name]} \else \label{#1} \fi #2&\lnP{:=}#3,\\#4&\lnP{:=}#5,\\#6&\lnP{:=}#7,\\#8&\lnP{:=}#9,  \end{split}\end{equation} }
		\newcommand{\envHLineCEDefNm}[9][*]{ \begin{equation}\begin{split} \ifx*#1 \text{[ERROR;need label name]} \else \label{#1} \fi #2&\lnP{:=}#3\\#4&\lnP{:=}#5\\#6&\lnP{:=}#7\\#8&\lnP{:=}#9  \end{split}\end{equation} }
		\newcommand{\envHLineCEDefNmPd}[9][*]{ \begin{equation}\begin{split} \ifx*#1 \text{[ERROR;need label name]} \else \label{#1} \fi #2&\lnP{:=}#3,\\#4&\lnP{:=}#5,\\#6&\lnP{:=}#7,\\#8&\lnP{:=}#9.  \end{split}\end{equation} }
	\newcommand{\envMLineTPd}[3][*]{ \ifx*#1 \begin{multline*} #2\lnP{=}#3.\end{multline*} \else \begin{multline} \label{#1} #2\lnP{=}#3.\end{multline} \fi }
	\newcommand{\envMLineTCm}[3][*]{ \ifx*#1 \begin{multline*} #2\lnP{=}#3,\end{multline*} \else \begin{multline} \label{#1} #2\lnP{=}#3,\end{multline} \fi }
		\newcommand{\envMLinePd}{\envMLineTPd}
		\newcommand{\envMLineCm}{\envMLineTCm}
		
		\newcommand{\envHLineTCmPt}[4][*]{\glcondEnvLineHead{#1} #3&#2&#4\glcondEnvLineTailCm{#1}}

			\newcommand{\envHLineCmPt}{\envHLineTCmPt}

		\newcommand{\envHLineFPdPte}[8][*]{\glcondEnvLineHead{#1} #2&#3&#4\\&#5&#6\nonumber \\&#7&#8\nonumber\glcondEnvLineTailPd{#1}}
			\newcommand{\lccondPar}[1]{\ifx#1p \\ \fi}
			\newcommand{\HLineCECm}[9][?]{#2&=&#3,\nonumber\\#4&=&#5,\nonumber\\#6&=&#7,\nonumber\\#8&=&#9,\nonumber\lccondPar{#1}}
			\newcommand{\HLineCEPd}[9][?]{#2&=&#3,\nonumber\\#4&=&#5,\nonumber\\#6&=&#7,\nonumber\\#8&=&#9.\nonumber\lccondPar{#1}}

			\newcommand{\OTLineCThCm}[4][?]{#2&=&#3=#4,\nonumber \ifx#1p \\ \fi}
			\newcommand{\OTLineCThPd}[4][?]{#2&=&#3=#4.\nonumber \ifx#1p \\ \fi}
			\newcommand{\OTLineCThCmDef}[4][?]{#2&:=&#3=#4,\nonumber \ifx#1p \\ \fi}
			\newcommand{\OTLineCThDefPd}[4][?]{#2&:=&#3=#4.\nonumber \ifx#1p \\ \fi}
			\newcommand{\OTLineCSCm}[7][?]{#2&=&#3=#4\nonumber#5&=&#6=#7,\nonumber \ifx#1p \\ \fi}
			\newcommand{\OTLineCSPd}[7][?]{#2&=&#3=#4\nonumber#5&=&#6=#7.\nonumber \ifx#1p \\ \fi}
	\newcommand{\envMLineTCmPt}[4][*]{ \ifx*#1 \begin{multline*} #3\lnP{#2}#4,\end{multline*} \else \begin{multline} \label{#1} #3\lnP{#2}#4,\end{multline} \fi }
	\newcommand{\envMLineTPdPt}[4][*]{ \ifx*#1 \begin{multline*} #3\lnP{#2}#4.\end{multline*} \else \begin{multline} \label{#1} #3\lnP{#2}#4.\end{multline} \fi }

\theoremstyle{plain}
\newtheorem{theorem}{THEOREM}[section]
\newtheorem{proposition}[theorem]{PROPOSITION}
\newtheorem{lemma}[theorem]{LEMMA}
\newtheorem{corollary}[theorem]{COROLLARY}
\theoremstyle{definition}

\theoremstyle{remark}
\newtheorem{remark}[theorem]{REMARK}
\allowdisplaybreaks[4]
\numberwithin{equation}{section}

	\newcommand{\lccondBibitem}[3][]{ \if ?#2 \bibitem{#3} \else \bibitem[#2]{#3} \fi}
	\newcommand{\refPaper}[8][?]{
			\lccondBibitem{#1}{#2}%bibitem		\bibitem[#1]{#2}  		%bibitem
				#3,			%author
				\emph{#4}, 	%artName
				#5\ 			%jouName 
				{\bf #6},		%volume
				#7,			%year
				#8.			%page
		}
	
	\newcommand{\refPaperRep}[9][?]{
			\lccondBibitem{#1}{#2}%bibitem
				#3,			%author
				\emph{#4}, 	%artName
				#5\ 			%jouName 
				{\bf #6},		%volume
				#7,			%year
				#8			%page
				; reprinted in #9	%reprinted information
		}
	
	\newcommand{\refBook}[7][?]{
			\lccondBibitem{#1}{#2}%bibitem
				#3,			%author
				\emph{#4}, 	%bookName
				#5,			%publisher 
				#6,			%pubPlace
				#7.			%pubYear
		}
	\newcommand{\refPaperAlm}[5][?]{
			\lccondBibitem{#1}{#2}%bibitem
				#3,	 		%author
				\emph{#4}, 	%artName
				#5		%etc 
		}

	\newcommand{\pcstSpForRefThm}{\ }
	\newcommand{\refThm}[2][?]{\ifx?#1\refH{Theorem\pcstSpForRefThm}{#2}\else\ifx s#1\refH{Theorems\pcstSpForRefThm}{#2}\else{[argument error]}\fi\fi}
	\newcommand{\refProp}[2][?]{\ifx?#1\refH{Proposition\pcstSpForRefThm}{#2}\else\ifx s#1\refH{Propositions\pcstSpForRefThm}{#2}\else{[argument error]}\fi\fi}
	\newcommand{\refLem}[2][?]{\ifx?#1\refH{Lemma\pcstSpForRefThm}{#2}\else\ifx s#1\refH{Lemmas\pcstSpForRefThm}{#2}\else{[argument error]}\fi\fi}
	\newcommand{\refCor}[2][?]{\ifx?#1\refH{Corollary\pcstSpForRefThm}{#2}\else\ifx s#1\refH{Corollaries\pcstSpForRefThm}{#2}\else{[argument error]}\fi\fi}
	\newcommand{\refDef}[2][?]{\ifx?#1\refH{Definition\pcstSpForRefThm}{#2}\else\ifx s#1\refH{Definitions\pcstSpForRefThm}{#2}\else{[argument error]}\fi\fi}
	\newcommand{\refRem}[2][?]{\ifx?#1\refH{Remark\pcstSpForRefThm}{#2}\else\ifx s#1\refH{Remarks\pcstSpForRefThm}{#2}\else{[argument error]}\fi\fi}
	\newcommand{\refTab}[2][?]{\ifx?#1\refH{Table\pcstSpForRefThm}{#2}\else\ifx s#1\refH{Tables\pcstSpForRefThm}{#2}\else{[argument error]}\fi\fi}

	\newcommand{\refPropF}[1]{\refF{Proposition}{#1}}

	\DeclareFontEncoding{OT2}{}{}
	\DeclareFontSubstitution{OT2}{cmr}{m}{n}
	\DeclareFontFamily{OT2}{cmr}{\hyphenchar\font45}
	\DeclareFontShape{OT2}{cmr}{m}{n}{<5><6><7><8><9>gen*wncyr <10><10.95><12><14.4><17.28><20.74><24.88>wncyr10}{}
	\DeclareFontShape{OT2}{cmr}{b}{n}{<5><6><7><8><9>gen*wncyb<10><10.95><12><14.4><17.28><20.74><24.88>wncyb10}{}
	\DeclareMathAlphabet{\mathcyr}{OT2}{cmr}{m}{n}
	\DeclareMathAlphabet{\mathcyb}{OT2}{cmr}{b}{n}
	\SetMathAlphabet{\mathcyr}{bold}{OT2}{cmr}{b}{n}

	\newcommand{\sh}{\mathcyr{sh}}

	\newcommand{\letGfcDZV}{\mathfrak{D}} %for using in the text

	\newcommand{\gfcDZV}[3][n]{\idFc[#1]{\letGfcDZV}{#2}{#3}}
	\newcommand{\gfcDPL}[4][n]{\idFc[#1]{\letGfcDZV\mathfrak{L}}{#2}{#3;#4}}
	\newcommand{\gfcTZV}[3][n]{\idFc[#1]{\mathfrak{T}}{#2}{#3}}

	\newcommand{\gfcTPL}[4][n]{\idFc[#1]{\mathfrak{TL}}{#2}{#3;#4}}

	\newcommand{\lue}{e}
	\newcommand{\lsetS}[2][?]{\id{U}{#2}}
	\newcommand{\lfcCst}[2][n]{\Fc[#1]{C_0}{#2}}
	\newcommand{\lmgfcDZV}[3][n]{\idFc[#1]{ \mathfrak{D}^{\sh} }{#2}{#3}}
	\newcommand{\lmgfcTZV}[3][n]{\idFc[#1]{ \mathfrak{T}^{\sh} }{#2}{#3}}
	\newcommand{\lfcSPLi}[5][n]{\ipFc[#1]{\mathfrak{PL}}{#2}{(#3)}{#4;#5}}
	\newcommand{\lfcSPZ}[4][n]{\ipFc[#1]{\mathfrak{P}}{#2}{(#3)}{#4}}
	\newcommand{\lmfcZ}[2][n]{\Fc[#1]{\zeta^{\sh}}{#2}}
	\newcommand{\lfcST}[2][n]{\Fc[#1]{\mathfrak{S}}{#2}}
	\newcommand{\lfcSTid}[3][n]{\idFc[#1]{\mathfrak{s}}{#2}{#3}}
	\newcommand{\lrelEDS}[4][?]{\ipFc{r}{#2}{(#3)}{#4}}
	\newcommand{\lclsEDS}[3][?]{\ip{R}{#2}{(#3)}}
	\newcommand{\loclsEDS}[3][?]{\ip{R}{#2}{#3}}
	\newcommand{\lalgFZV}[1][?]{ \ifx ?#1 \mathfrak{H} \else \mathfrak{H}^{#1} \fi }
	\newcommand{\lalgFZVs}[1][?]{ \ifx ?#1 \mathfrak{H}_{\sh} \else \mathfrak{H}_{\sh}^{#1} \fi }
	\newcommand{\lmpZV}[2][n]{\Fc[#1]{Z}{#2}}
	\newcommand{\lmpZVs}[2][n]{\pwFc[#1]{Z}{\sh}{#2}}
	\newcommand{\lmpZVi}[3][n]{\ipFc[#1]{Z}{#2}{\sh}{#3}}
	\newcommand{\lmpReg}[2][n]{\idFc[#1]{\mathrm{reg}}{\sh}{#2}}
	\newcommand{\lmpRegT}[2][n]{\ipFc[#1]{\mathrm{reg}}{\sh}{T}{#2}}
	\newcommand{\lSm}[2][?]{\ifx t#1 \sum_{#2}^{\prime} \else \sideset{}{^\prime}\sum_{#2} \fi}
	\newcommand{\lsh}{\ \sh\ }
	
	\newcommand{\lrpTx}[1]{(resp.\;#1)}
	\newcommand{\rTx}[2][?]{ \raise1ex\hbox{$ \displaystyle#2$} }
	
\allowdisplaybreaks[4]
	\title{\mainTitle}
	\author{\authorName}
	\date{}

\begin{document}
%--maketitle
\maketitle
%--abstract
\begin{abstract}
Extended double shuffle relations for multiple zeta values are 
	obtained by using the fact that
	any product of regularized multiple zeta values has two different representations.
In this paper,
	we give two formulas for the generating function of the triple zeta values of any fixed weight
	by the use of the extended double shuffle relations
	obtained as two-fold products of double and single zeta values and also as three-fold products of single zeta values.
As applications of the formulas,
	we also obtain parameterized, weighted, and restricted sum formulas for triple zeta values.

%\mbox{}\\ \\
%{\bf Key Words and Phrases:}
%triple zeta value, multiple zeta value, multiple polylogarithm, 
%extended double shuffle relation, sum formula
%
%\mbox{}\\
%{\bf MSC classes:}
%11M32
\end{abstract}

%--body
\section{\sectOne} \label{sectOne}
A \emph{multiple zeta value}, 
	sometimes called a \emph{multiple harmonic series} or a \emph{Euler sum},
	is defined by the infinite series
	\envHLineDef[1_Def_MZV]
	{
		\fcZ{l_1,l_2,\ldots,l_n}
	}
	{
		\Sm{m_1>m_2>\cdots>m_n>0} \opF{1}{ \pw{m_1}{l_1}\pw{m_2}{l_2}\cdots\pw{m_n}{l_n} }
	}
	for any multi-index $(l_1,l_2,\ldots,l_n)$ of positive integers satisfying the condition $l_1\geq2$, 
	which is necessary for convergence.
The integers $l=l_1+\cdots+l_n$ and $n$ are called the weight and depth of the multiple zeta value, respectively.
\emph{Regularized multiple zeta values},
	which were introduced in \cite{IKZ06} in order to deal with the divergent series as meaningful objects,
	consist of multiple zeta values 
	and
	certain real numbers which are obtained by adjoining 
	infinite values $\fcZ{1,l_2,\ldots,l_n}$ to suitable combinations of multiple zeta values.
It is known that there are two large collections of linear relations over the rational numbers $\setQ$ among multiple zeta values:
	Kawashima's relations \cite{Kawashima09} and \emph{extended double shuffle relations} \cite{IKZ06}.
The latter relations, 
	which play an important role in the present paper,
	are derived from the two sets of relations
	which are obtained by using the fact that any product of regularized multiple zeta values 
	has two different $\setZ$-linear combinations.
The fact is due to
	two expressions of multiple zeta values,
	which are formed by integration and summation.
We refer to
	the relations induced from the integral \lrpTx{summation} expression \emph{extended shuffle} \lrpTx{\emph{harmonic}} \emph{relations}.
(See \cite[(1.2)]{IKZ06} for the integral expression, and \cite[(1.1)]{IKZ06} or \refEq{1_Def_MZV} for the summation expression.)
The extended shuffle relations
	are also obtained as either the integral expression or as the partial fraction decomposition \cite{KMT11}.
It should be noted 
	that there are two kinds of regularized multiple zeta value: one kind related to shuffle and one kind related to harmonic, 
	and that there are several equivalent definitions of extended double shuffle relations (see \cite[Theorem 2]{IKZ06}).
The kind of regularized multiple zeta value and the definition of extended double shuffle relations which are used in the present paper 
	concern the shuffle regularization and \cite[Theorem 2(iv)]{IKZ06}, respectively.

Double zeta values, already studied by Euler \cite{Euler1775}, 
	are the multiple zeta values of depth $2$.
Recently, Gangl, Kaneko, and Zagier \cite{GKZ06} showed that
	extended double shuffle relations for double zeta values yield an elegant formula for 
	the generating function of the double zeta values of any weight $l$.
The aims of the present paper are to generalize their formula to triple zeta values and to give some applications of the generalized obtained formulas. 
Therefore, we will first review this prior work.
	
We define the generating function of the double zeta values of weight $l$ as
	\envHLineDefPd[1_Plane_DefGfcDZV]
	{
		\gfcDZV{l}{x,y} 
	}
	{
		\tpSm{l_1\geq2,l_2\geq1}{l_1+l_2=l} x^{l_1-1}y^{l_2-1} \fcZ{l_1,l_2}
	}
(Gangl, Kaneko, and Zagier \cite{GKZ06} use $\mathfrak{Z}$ rather than $\letGfcDZV$ to denote the generating function;
	however, we adopt the symbol $\mathfrak{D}$ to emphasize that $\gfcDZV{l}{x,y}$ is related to double zeta values.)
They
	formally derived the formula \cite[(26)]{GKZ06} from the relations \cite[(22)]{GKZ06}.
Through 
	the $\setR$-realization \cite[(24)]{GKZ06} with $\kappa=-\bkR[b]{\fcZ{k-1,1}+\SmT{j=2}{k-1}\fcZ{j,k-j}}=-\bkR[b]{\fcZ{k-1,1}+\fcZ{k}}$
	where the last equality follows from \refEq{1_Plane_EqSFofDZV} below,
	\cite[(22)]{GKZ06} become the extended double shuffle relations for double zeta values.
Thus these relations yield
	\envHLinePd[1_Eq_gfcDZ]
	{
		\gfcDZV{l}{x+y,y} + \gfcDZV{l}{y+x,x} 
	}
	{
		\gfcDZV{l}{x,y} + \gfcDZV{l}{y,x} + \opF{x^{l-1}-y^{l-1}}{x-y}\fcZ{l}
	}
The formula \refEq{1_Eq_gfcDZ} can be considered to be a \emph{parameterized sum formula} or 
	a parameterized analogue of the following sum formula given by Euler \cite{Euler1775},
	\envOTLineCm[1_Plane_EqSFofDZV]
	{
		\tpSm{l_1\geq2,l_2\geq1}{l_1+l_2=l} \fcZ{l_1,l_2}
	}
	{
		\gfcDZV{l}{1,1} 
	}
	{
		\fcZ{l}
	}
	because \refEq{1_Eq_gfcDZ} yields some known sum formulas if the appropriate values are substituted for the parameters $x$ and $y$.
In fact, \refEq{1_Eq_gfcDZ} with $(x,y)=(1,0)$ and $(1,1)$ respectively give the original sum formula \refEq{1_Plane_EqSFofDZV}
	and the \emph{weighted sum formula} \cite[(10)]{OZu08}.
If $l$ is even, 
	\refEq{1_Eq_gfcDZ} with  $(x,y)=(-1,1)$ yields 
	\envM
	{
		\gfcDZV{l}{-1,1}
	}
	{
		-\opF[s]{\fcZ{l}}{2}
	}
	because 
	\envMCm{ \gfcDZV{l}{-x, -y} }{ \gfcDZV{l}{x, y} }
	and 
	we obtain the \emph{restricted sum formulas} \cite[(4)]{GKZ06}
	\envHLineCFPd[1_Plane_EqRSFofDZV]
	{
		\tpSm{l_1\geq2,l_2\geq1}{l_1+l_2=l \atop l_1,l_2:even} \fcZ{l_1,l_2}
	}
	{
		\opF{\gfcDZV{l}{1,1} - \gfcDZV{l}{-1,1}}{2}
	\lnP{=}
		\opF{3}{4} \fcZ{l}
	}
	{
		\tpSm{l_1\geq2,l_2\geq1}{l_1+l_2=l \atop l_1,l_2:odd} \fcZ{l_1,l_2}
	}
	{
		\opF{\gfcDZV{l}{1,1} + \gfcDZV{l}{-1,1}}{2}
	\lnP{=}
		\opF{1}{4} \fcZ{l}
	}
(See \cite{Nakamura09} for different types of restricted sum formula of double zeta values.)

Let $\gfcTZV{l}{x_1,x_2,x_3}$ denote the generating function of the triple zeta values of weight $l$, which is defined by
	\envHLineDefPd[1_Plane_DefGfcTZV]
	{
		\gfcTZV{l}{x_1,x_2,x_3}
	}
	{
		\tpSm{l_1\geq2,l_2,l_3\geq1}{l_1+l_2+l_3=l} x_1^{l_1-1}x_2^{l_2-1}x_3^{l_3-1} \fcZ{l_1,l_2,l_3}
	}
The first purpose of this paper is 
	to give two formulas for $\gfcTZV{l}{x_1,x_2,x_3}$
	by making use of the two classes of extended double relations of the triple zeta values of weight $l$.	
One of the classes, $\lclsEDS{l}{2,1}$, is 
	the set of relations $\lrelEDS{l}{2,1}{p,q,r}$ obtained as two-fold products $\fcZ{p,q}\fcZ{r}$,
	and
	the another class, $\lclsEDS{l}{1,1,1}$, is the set of relations $\lrelEDS{l}{1,1,1}{p,q,r}$ obtained as three-fold products $\fcZ{p}\fcZ{q}\fcZ{r}$,
	where $p, q, r$ are positive integers satisfying $l=p+q+r$
	and divergent series are considered to be regularized multiple zeta values.		
(See \refEq{2.2_Def_CtDTZV} for explicit expressions of these values
	and also the coefficients of $x_1^{p-1}x_2^{q-1}x_3^{r-1}$ of \refEq{3_Prop1i_Eq1} and \refEq{3_Prop1ii_Eq1}
	for those of $\lrelEDS{l}{2,1}{p,q,r}$ and  $\lrelEDS{l}{1,1,1}{p,q,r}$, respectively.)	
Only the relations
	\envMDef{ \loclsEDS{l}{3} }{ \bkB[b]{ \lrelEDS{l}{2,1}{p,q,r}\in\lclsEDS{l}{2,1} \mVert (p,q,r) \neq (1,l-2,1) } }
	are treated as extended double shuffle relations for triple zeta values in \cite[Theorem 2(iv)]{IKZ06},
	but we also regard $\lrelEDS{l}{2,1}{1,l-2,1}$ and $\lclsEDS{l}{1,1,1}$ as extended double shuffle relations in the present paper
	since 
	$\lrelEDS{l}{2,1}{1,l-2,1}$ and $\lclsEDS{l}{1,1,1}$ can be obtained in the same manner as the original relations $\loclsEDS{l}{3}$.
	
We prepare by defining some notation before stating the two formulas.	
Let $i,j,k$ be distinct integers satisfying $1\leq i,j,k \leq3$.
For variables $x_1,x_2,x_3$, 
	we define
	\envPLinePd[1_Def_ijk]
	{
		x_{ij} \lnP{:=} x_i+x_j,	\qquad
		x_{ijk} \lnP{:=} x_i+x_j+x_k	\ \bkR[n]{=x_1+x_2+x_3 } 
	}
Let $\gpSym{3}$ be the symmetric group of degree $3$ and $\lue$ be its identity element.
We respectively denote the cycle permutations 
	\envMO
	{
		\envSMatT{i&j&k}{j&i&k}
	}
	and
	\envMO
	{
		\envSMatT{i&j&k}{j&k&i}
	}
	by  $(ij)$ and $(ijk)$,
	and the subsets $\bkB{\lue,(ij)}$ and $\bkB{\lue,(ij),(ijk)}$ of $\gpSym{3}$
	by $\lsetS{ij}$ and $\lsetS{ijk}$.
The alternating group $\bkB{\lue, (123), (132)}$ of degree $3$ is denoted  by $\gpAlt{3}$.	

The previously mentioned two formulas are the following equations,
	 \refEq{1_Thm1i_Eq1} and \refEq{1_Thm1ii_Eq1},
	which shall be proved by the use of the two classes 
	$\lclsEDS{l}{2,1}$ and $\lclsEDS{l}{1,1,1}$, respectively.
	
%A
%:1_Thm1
\begin{theorem}\label{1_Thm1}
Let $l$ be an integer such that $l\geq4$, and $x_1,x_2,x_3$ be variables. Then the following two equations hold:
\mbox{}\\{\bf (i)}
	\envLinePd[1_Thm1i_Eq1]
	{
		\gfcTZV{l}{x_{13},x_{23},x_{3}}+\gfcTZV{l}{x_{13},x_{32},x_{2}}+\gfcTZV{l}{x_{31},x_{1},x_{2}}
		+
		\opF{\gfcDZV{l}{x_{3},x_{2}}}{x_{3}} 
	}
	{
		\Sm{\sig\in\lsetS{321}} \gfcTZV{l}{x_{\sig(1)},x_{\sig(2)},x_{\sig(3)}}
		+
		\gfcTZV{l}{x_{3},x_{23},x_{3}}+\gfcTZV{l}{x_{3},x_{32},x_{2}} 
		\lnAHs{140}
		+
		\Sm{\sig\in\lsetS{31}} \opF{\gfcDZV{l}{x_{\sig(1)},x_{2}}}{x_{\sig(1)} -  x_{\sig(3)}} 
		+
		\Sm{\sig\in\lsetS{32}} \opF{\gfcDZV{l}{x_{1},x_{\sig(2)}}}{x_{\sig(2)} - x_{\sig(3)}} 
	}
{\bf (ii)} 
	\envLinePd[1_Thm1ii_Eq1]
	{
		\Sm{\sig\in\gpSym{3}} \gfcTZV{l}{x_{\sig(1)\sig(2)\sig(3)},x_{\sig(2)\sig(3)},x_{\sig(3)}}
		+
		\Sm{\tau\in\gpAlt{3}} \opF{\gfcDZV{l}{x_{\tau(3)},x_{\tau(2)}}}{x_{\tau(3)}} 
	}
	{
		\Sm{\sig\in\gpSym{3}}
		\bkS[g]{
			\gfcTZV{l}{x_{\sig(1)},x_{\sig(2)},x_{\sig(3)}}
			+
			\opF{\gfcDZV{l}{x_{\sig(1)},x_{\sig(2)}}}{x_{\sig(1)} - \pw{x_{\sig(3)}}{1}} 
			+
			\opF{\gfcDZV{l}{x_{\sig(1)},x_{\sig(2)}}}{x_{\sig(2)} - \pw{x_{\sig(3)}}{1}} 
		}
		\lnAHs{-5}
		+
		\Sm{\tau\in\gpAlt{3}} 
		\bkS[g]{
			\gfcTZV{l}{x_{\tau(3)},x_{\tau(2)\tau(3)},x_{\tau(3)}}+\gfcTZV{l}{x_{\tau(3)},x_{\tau(3)\tau(2)},x_{\tau(2)}}		
			+
			\bkR[B]{\ip{x}{\tau(1)}{l-1} \PdT{i=2}{3} \opF{1}{\pw{x_{\tau(1)}-x_{\tau(i)}}{1}}} \fcZ{l}
		}
	}
\end{theorem}
%Z

The formulas \refEq{1_Thm1i_Eq1} and \refEq{1_Thm1ii_Eq1} include many double zeta values
	and so
	it is difficult to claim that \refEq{1_Thm1i_Eq1} and \refEq{1_Thm1ii_Eq1} are analogues of the following sum formula 
	for triple zeta values which was proved in \cite{HM96}
	(see \cite{Granville97,Zagier95} for the case of multiple zeta values of any depth):
	\envOTLinePd[1_Plane_EqSFofTZV]
	{
		\tpSm{l_1\geq2,l_2,l_3\geq1}{l_1+l_2+l_3=l} \fcZ{l_1,l_2,l_3}
	}
	{
		\gfcTZV{l}{1,1,1}
	}
	{
		\fcZ{l}
	} 
	
The second purpose of the present paper is to derive various formulas which are analogues of the sum formula \refEq{1_Plane_EqSFofTZV}
	 from \refThm{1_Thm1}.
One of the formulas, stated in \refThm{1_Thm2} below, is a parameterized sum formula.
It does not include double zeta values,
	and it becomes the original sum formula \refEq{1_Plane_EqSFofTZV} if $(1,0,0)$ is substituted for $(x_1,x_2,x_3)$.
Furthermore, the parameterized sum formula below yields weighted sum formulas which contain  the result of Guo and Xie \cite[Theorem 1.1]{GX09} 
	in the case of triple zeta values (see \refCor{4_Cor1}).
%A
%:1_Thm2
\begin{theorem}\label{1_Thm2}
Let $l$ be an integer such that $l\geq4$, and $x_1,x_2,x_3$ be variables.
We have
	\envLinePd[1_Thm2_Eq1]
	{
		\Sm{\sig\in\gpSym{3}} \gfcTZV{l}{x_{\sig(1)\sig(2)\sig(3)},x_{\sig(2)\sig(3)},x_{\sig(3)}}
	}
	{
		\Sm{\tau\in\gpAlt{3}} 
		\bkS[G]{
			\gfcTZV{l}{ x_{\tau(1)\tau(3)},x_{\tau(2)\tau(3)},x_{\tau(3)}}
			+
			\gfcTZV{l}{ x_{\tau(1)\tau(3)},x_{\tau(3)\tau(2)},x_{\tau(2)}}
			\lnAHs{25}
			+
			\gfcTZV{l}{x_{\tau(3)\tau(1)},x_{\tau(1)},x_{\tau(2)}}
			-
			\gfcTZV{l}{x_{\tau(1)},x_{\tau(2)},x_{\tau(3)}}
			+
			\bkR[g]{\ip{x}{\tau(1)}{l-1} \nPdT{i=2}{3} \opF{1}{\pw{x_{\tau(1)}-x_{\tau(i)}}{1}}} \fcZ{l}
		}
	}
In other words, 
	\envLinePd[1_Thm2_Eq2]
	{
		\tpSm{l_1\geq2,l_2,l_3\geq1}{l_1+l_2+l_3=l}  
		\bkS[G]{ 
			\Sm{\sig\in\gpSym{3}} x_{\sig(1)\sig(2)\sig(3)}^{l_1-1}x_{\sig(2)\sig(3)}^{l_2-1}x_{\sig(3)}^{l_3-1} 
			+
			\Sm{\tau\in\gpAlt{3}} \bkR[g]{
				x_{\tau(1)}^{l_1-1}x_{\tau(2)}^{l_2-1}x_{\tau(3)}^{l_3-1}
				\lnAHs[]{10}
				-
				x_{\tau(1)\tau(3)}^{l_1-1}x_{\tau(2)\tau(3)}^{l_2-1}x_{\tau(3)}^{l_3-1}
				-
				x_{\tau(1)\tau(3)}^{l_1-1}x_{\tau(3)\tau(2)}^{l_2-1}x_{\tau(2)}^{l_3-1}
				-
				x_{\tau(3)\tau(1)}^{l_1-1}x_{\tau(1)}^{l_2-1}x_{\tau(2)}^{l_3-1}
			}
		} 
		\fcZ{l_1,l_2,l_3}		\nonumber
	}
	{
		\bkR[G]{ \tpSm{l_1,l_2,l_3\geq1}{l_1+l_2+l_3=l} x_1^{l_1-1}x_2^{l_2-1}x_3^{l_3-1} } \fcZ{l}
	}
\end{theorem}
%Z

The rest of the formulas which are of interest that
	can be obtained using \refThm{1_Thm1} (and \refThm{1_Thm2}) 
	are restricted analogues of the sum formula \refEq{1_Plane_EqSFofTZV} (see \refEq{1_Thm3i_Eq1} and \refEq{1_Thm3ii_Eq1} below).
The first equation of \refEq{1_Thm3i_Eq1}, 
	which is equivalent to the second one by \refEq{1_Plane_EqSFofTZV}, 
	is the result of Shen and Cai \cite[Theorem 1]{SC12}.
We also calculate restricted analogues of the formulas given by Granville, Hoffman, and Ohno in \refSect{sectFive} (see \refProp{5_Prop1}).
	
%A
%:1_Thm3
\begin{theorem}[\text{cf. \cite[Theorem 1]{SC12}}]\label{1_Thm3}
Let $l$ be an integer such that $l\geq4$.
For any condition $P(l_1,\ldots,l_n)$ of positive integers $l_1,\ldots,l_n$,
	we mean by $\pSm[t]{P(l_1,\ldots,l_n)}$ summing over all integers 
	satisfying $l_1\geq2$, $l_2,\ldots,l_n\geq1,$ $l=l_1+\cdots+l_n$, and $P(l_1,\ldots,l_n)$.
\mbox{}\\{\bf (i)} If $l$ is even, then
	\envHLineCSPd[1_Thm3i_Eq1]
	{
		\pSm{l_1,l_2,l_3:even} \fcZ{l_1,l_2,l_3}
	}
	{
		\opF{5}{8} \fcZ{l} -\opF{1}{4}\fcZ{l-2}\fcZ{2}
	}
	{
		\bkR{ \rTx{ \pSm{l_1:even\atop l_2,l_3:odd} + \pSm{l_2:even\atop l_1,l_3:odd} + \pSm{l_3:even\atop l_1,l_2:odd} } } \fcZ{l_1,l_2,l_3}
	}
	{
		\opF{3}{8} \fcZ{l} + \opF{1}{4}\fcZ{l-2}\fcZ{2}
	}
	{
		\bkR{ \rTx{\pSm{l_1:even\atop l_2,l_3:odd} - \pSm{l_3:even\atop l_1,l_2:odd}}  } \fcZ{l_1,l_2,l_3}
	}
	{
		-\opF{1}{4} \fcZ{l} + \opF{1}{2}\fcZ{l-2}\fcZ{2}
	}	
{\bf (ii)} If $l$ is odd, then
	\envHLineCSPd[1_Thm3ii_Eq1]
	{
		\pSm{l_1,l_2,l_3:odd} \fcZ{l_1,l_2,l_3} +  \opF{1}{2} \pSm{l_2:even \atop l_1:odd}\fcZ{l_1,l_2}
	}
	{
		\opF{3}{8}\fcZ{l} - \opF{1}{4}\fcZ{2,l-2} 
	}
	{
		\pSm{l_1,l_3:even \atop l_2:odd} \fcZ{l_1,l_2,l_3} + \opF{1}{2} \pSm{l_2:even \atop l_1:odd}\fcZ{l_1,l_2}
	}
	{
		\opF{5}{8}\fcZ{l} + \opF{1}{4}\fcZ{2,l-2} 
	}
	{
		\bkR{ \rTx{\pSm{l_1,l_2:even\atop l_3:odd} +  \pSm{l_2,l_3:even\atop l_1:odd}}  } \fcZ{l_1,l_2,l_3}
		+
		\pSm{l_1:even \atop l_2:odd}\fcZ{l_1,l_2}
	}
	{
		\fcZ{l}
	}
\end{theorem}
%Z
We will first outline how the theorems can be proved. 
Let $\lgP{l_1,\ldots,l_n}{z_1,\ldots,z_n}$ be the \emph{multiple polylogarithm} defined by
	\envHLineDef[1_Def_lgP]
	{
		\lgP{l_1,\ldots,l_n}{z_1,\ldots,z_n}
	}
	{
		\Sm{m_1>\cdots>m_n>0} \opF{ \pw{z_1}{m_1-m_2} \cdots \pw{z_{n-1}}{m_{n-1}-m_n}\pw{z_n}{m_n} }{ \pw{m_1}{l_1}\cdots\pw{m_{n-1}}{l_{n-1}}\pw{m_n}{l_n} }
	}
	for any multi-index $(l_1,\ldots,l_n)$ of positive integers and complex numbers $z_i\ (i=1,\ldots,n)$ such that $\abs{z_i}<1$.
Unlike multiple zeta values,
	multiple polylogarithms with $l_1=1$ converge,  
	because $\abs{z_i}<1$.	
We define the generating functions of the double and triple polylogarithms of weight $l$ as
	\envHLineCFDefPd[1_Plane_DefGfcPL]
	{
		\gfcDPL{l}{x_1,x_2}{z_1,z_2}
	}
	{
		\tpSm{l_1,l_2\geq1}{l_1+l_2=l} x_1^{l_1-1}x_2^{l_2-1} \lgP{l_1,l_2}{z_1,z_2}
	}
	{
		\gfcTPL{l}{x_1,x_2,x_3}{z_1,z_2,z_3}
	}
	{
		\tpSm{l_1,l_2,l_3\geq1}{l_1+l_2+l_3=l} x_1^{l_1-1}x_2^{l_2-1}x_3^{l_3-1} \lgP{l_1,l_2,l_3}{z_1,z_2,z_3}
	}
Note that the summations in \refEq{1_Plane_DefGfcPL} allow $l_1=1$, unlike those in \refEq{1_Plane_DefGfcDZV} and \refEq{1_Plane_DefGfcTZV}. 
In order to prove the theorems,
	we first give formulas for some functions expressed in terms of \refEq{1_Plane_DefGfcPL}
	with $z_i\in\bkB{z,z^2,z^3}$ in \refProp{2.1_Prop1}.
The formulas relate to shuffle and harmonic relations (see \refRem{2.1_Rem1} for details).
In \refProp{2.2_Prop1}, we also give asymptotic properties of the functions which appear in the formulas of \refProp{2.1_Prop1}.
Next we calculate limits of the formulas of \refProp{2.1_Prop1} as $z\nearrow1$ by the use of \refProp{2.2_Prop1}
	and obtain identities for the functions $\ipFc{F}{2}{\sh}{x_1,x_2} = \lmgfcDZV{l}{x_1,x_2}$ and $\ipFc{F}{3}{\sh}{x_1,x_2,x_3} =\lmgfcTZV{l}{x_1,x_2,x_3}$ in \refProp{3_Prop1}.
Here $\ipFc{F}{2}{\sh}{x_1,x_2}$ and $\ipFc{F}{3}{\sh}{x_1,x_2,x_3}$ are as defined in \cite[\S8]{IKZ06},
	and are the generating functions of the regularized double and triple zeta values, respectively 
	(see also \cite{IO08}, in which the generating function of triple zeta values are studied).
These generating functions are modified versions of $\gfcDZV{l}{x_1,x_2}$  and $\gfcTZV{l}{x_1,x_2,x_3}$, as we shall see in \refEq{2.2_Plane_DefMGFC},
	and so for simplicity 
	we respectively use $\lmgfcDZV{l}{x_1,x_2}$ and $\lmgfcTZV{l}{x_1,x_2,x_3}$ 
	instead of $\ipFc{F}{2}{\sh}{x_1,x_2}$ and $\ipFc{F}{3}{\sh}{x_1,x_2,x_3}$.
We then prove \refThm{1_Thm1} using \refProp{3_Prop1}
	and derive the other theorems from \refThm{1_Thm1}.
It is worth noting that
	Borwein and Girgensohn \cite{BG96} adopted a similar but not identical proof of a parity result regarding triple zeta values.
Instead of multiple polylogarithms $\lgP{l_1,\ldots,l_n}{z_1,\ldots,z_n}$, 
	they used the partial zeta sums 
	\envMDef{ \idFc{\zeta}{N}{l_1,\ldots, l_n} }{ \Sm[l]{N>m_1>\cdots>m_n>0}\opF[s]{1}{m_1^{l_1}\cdots m_n^{l_n}} }
	and considered the asymptotic properties of these sums as $N\to\infty$ (see their paper for details).
	
The present paper is organized as follows.
We respectively verify \refProp[s]{2.1_Prop1}, \ref{2.2_Prop1}, and \ref{3_Prop1} 
	in \refSect{sSectTwoO}, \refSect{sSectTwoT}, and the first half of \refSect{sectThree},
	and then 
	prove \refThm[s]{1_Thm1}, \ref{1_Thm2}, and \ref{1_Thm3} in the latter half of \refSect{sectThree}, \refSect{sectFour}, and \refSect{sectFive}.
We also give some weighted sum formulas in \refSect{sectFour}
	and restricted sum formulas in \refSect{sectFive} (see \refCor{4_Cor1} and \refProp{5_Prop1}, respectively).

\section{\sectTwo} \label{sectTwo}
\subsection{\sSectTwoO} \label{sSectTwoO}
In this subsection, 
	we give the formulas \refEq{2.1_Prop1i_Eq1} and \refEq{2.1_Prop1ii_Eq1} below
	for the generating functions \refEq{1_Plane_DefGfcPL} of the double and triple polylogarithms of any weight $l$.
To prove the formulas,
	we use the partial fraction expansion and the decomposition of summations
	which yield shuffle and harmonic relations, respectively;
	thus, the formulas correspond to these relations (see \refRem{2.1_Rem1} for details of the correspondences).
	
%A
%:2.1_Prop1
\begin{proposition}\label{2.1_Prop1}
Let $l$ be an integer such that $l\geq3$, 
	and $x_1,x_2,x_3$ be variables.
We define
	\envHLineCFCmDef
	{
		\lfcSPLi{l}{2}{x_1,x_2,x_3}{z_1,z_2}
	}
	{
		\tpSm{l_1\geq2, l_2 \geq 1}{l_1+l_2=l} \gfcDPL{l_1}{x_{1},x_{2}}{z_1,z_2}  x_3^{l_2-1} \lgP{l_2}{z_1}
	}
	{
		\lfcSPLi{l}{3}{x_1,x_2,x_3}{z}
	}
	{
		\tpSm[l]{l_1, l_2, l_3 \geq 1}{l_1+l_2+l_3=l} x_1^{l_1-1} x_2^{l_2-1} x_3^{l_3-1} \lgP{l_1}{z} \lgP{l_2}{z} \lgP{l_3}{z}
	}
	where $z, z_1,z_2$ are complex numbers such that $\abs{z},\abs{z_1},\abs{z_2}<1$.
\\{\bf (i)}
We have
	\envHLineCFPd[2.1_Prop1i_Eq1]
	{
		\lfcSPLi{l}{2}{x_{12},x_2,x_3}{z,z}
	}
	{
		\Sm{\sig\in\lsetS{321}} \gfcTPL{l}{x_{ \sig(1)\sig(2)\sig(3) },x_{ \sig(2)\sig(3) },x_{\sig(3)}}{z,z,z}
	}
	{
		\lfcSPLi{l}{2}{x_1,x_2,x_3}{z,z^2}	
	}
	{
		\Sm{\sig\in\lsetS{321}} \gfcTPL{l}{x_{\sig(1)},x_{\sig(2)},x_{\sig(3)}}{z,z^2,z^3}
		\lnAHs{0}
		+
		\Sm{\sig\in\lsetS{31}} \opF{\gfcDPL{l}{x_{\sig(1)},x_{2}}{z^2,z^3}}{x_{\sig(1)} -  x_{\sig(3)}} 
		+
		\Sm{\sig\in\lsetS{32}} \opF{\gfcDPL{l}{x_{1},x_{\sig(2)}}{z, z^3}}{x_{\sig(2)} -  x_{\sig(3)}} 	
	}
\\{\bf (ii)}
We have
	\envHLineThPd[2.1_Prop1ii_Eq1]
	{	
		\lfcSPLi{l}{3}{x_1,x_2,x_3}{z}
	}
	{
		\Sm{\sig\in\gpSym{3}} \gfcTPL{l}{x_{\sig(1)\sig(2)\sig(3)},x_{\sig(2)\sig(3)},x_{\sig(3)}}{z,z,z}
	}
	{
		\Sm{\sig\in\gpSym{3}}
		\bkS[G]{
			\gfcTPL{l}{x_{\sig(1)},x_{\sig(2)},x_{\sig(3)}}{z, z^2, z^3}
			\lnAHs{30}
			+
			\opF{ \gfcDPL{l}{x_{\sig(1)},x_{\sig(2)}}{z^2, z^3}}{x_{\sig(1)} - x_{\sig(3)}}
			+
			\opF{\gfcDPL{l}{x_{\sig(1)},x_{\sig(2)}}{z, z^3}}{x_{\sig(2)} - x_{\sig(3)}}
		}
		\lnAHs{0}
		+
		\Sm{\tau\in\gpAlt{3}} \bkR[g]{x_{\sig(1)}^{l-1} \PdT{i=2}{3} \opF{1}{x_{\sig(1)}-x_{\sig(i)}}} \lgP{l}{z^3}
	}
\end{proposition}
%Z
We remark that the left-hand sides of the two equations of \refEq{2.1_Prop1i_Eq1} are very similar, 
	but they are different in the argument of the function $\gfcDPL{l}{x_1,x_2}{z_1,z_2}$ appearing in the definition of $\lfcSPLi{l}{2}{x_1,x_2,x_3}{z_1,z_2}$:
in one the argument  is $(x_{12},x_{2};z,z)$ and in the other it is $(x_{1},x_{2};z,z^2)$.

We prepare a lemma to prove the proposition.
%A
%:2.1_Lem1
\begin{lemma}\label{2.1_Lem1}
For positive integers $k_1,k_2,k_3$, we have
	\envHLineCFNmePd
	{\label{2.1_Lem1Eq1}
		\lgP{k_1,k_2}{z,z^2}\lgP{k_3}{z}
	}
	{
		\Sm{\sig\in\lsetS{321}} \lgP{k_{\sig(1)},k_{\sig(2)},k_{\sig(3)}}{z,z^2,z^3} 
		\lnAHs[]{20}
		+ 
		\lgP{k_1+k_3,k_2}{z^2,z^3} + \lgP{k_1,k_2+k_3}{z,z^3}
		\nonumber
	}
	{\label{2.1_Lem1Eq2}
		\lgP{k_1}{z}\lgP{k_2}{z}\lgP{k_3}{z}
	}
	{
		\Sm{\sig\in\gpSym{3}} \lgP{k_{\sig(1)},k_{\sig(2)},k_{\sig(3)}}{z,z^2,z^3}
		\lnAHs[]{0}
		+
		\Sm{\tau\in\gpAlt{3}} \lgP{k_{\tau(1)}+k_{\tau(2)},k_{\tau(3)}}{z^2,z^3}
		+
		\Sm{\tau\in\gpAlt{3}} \lgP{k_{\tau(1)},k_{\tau(2)}+k_{\tau(3)}}{z,z^3}
		\lnAHs{20}
		+
		\lgP{k_1+k_2+k_3}{z^3}
		\nonumber
	}
\end{lemma}
%Z
%A
\envProof{
We see from \refEq{1_Def_lgP} that
	\envOTLinePd
	{
		\lgP{k_1,k_2}{z_1,z_2}\lgP{k_3}{z_3}
	}
	{
		\Sm{m_1>m_2>0} \opF{ \pw{z_1}{m_1-m_2}\pw{z_2}{m_2} }{ \pw{m_1}{k_1}\pw{m_2}{k_2} } \Sm{m_3>0} \opF{ \pw{z_3}{m_3} }{ \pw{m_3}{k_3} }
	}
	{
		\Sm{m_1>m_2>0 \atop m_3>0} \opF{ \pw{z_1}{m_1-m_2}\pw{z_2}{m_2}\pw{z_3}{m_3} }{ \pw{m_1}{k_1}\pw{m_2}{k_2} \pw{m_3}{k_3} } 
	}
Since we can decompose the summation $\Sm{m_1>m_2>0, m_3>0}$ as 
	\envHLineCm
	{
		\Sm{m_1>m_2>0 \atop m_3>0}
	}
	{
		\Sm{m_1>m_2> m_3>0} + \Sm{m_1>m_3> m_2>0} + \Sm{m_3>m_1> m_2>0} + \Sm{m_1=m_3> m_2>0} + \Sm{m_1>m_2= m_3>0} 
	}
	we obtain
	\envLinePd
	{
		\lgP{k_1,k_2}{z_1,z_2}\lgP{k_3}{z_3}
	}
	{
		\lgP{k_1,k_2,k_3}{z_1,z_2,z_2z_3} + \lgP{k_1,k_3,k_2}{z_1,z_1z_3,z_2z_3} + \lgP{k_3,k_1,k_2}{z_3,z_1z_3,z_2z_3} 
		\lnAHs{0}
		+ 
		\lgP{k_1+k_3,k_2}{z_1z_3,z_2z_3} + \lgP{k_1,k_2+k_3}{z_1,z_2z_3}
	}
Substituting $(z,z^2,z)$ for $(z_1,z_2,z_3)$ gives \refEq{2.1_Lem1Eq1}.
By the decomposition
	\envHLineCm
	{
		\Sm{m_1, m_2, m_3>0}
	}
	{
		\Sm{\sig\in\gpSym{3}} \Sm{m_{\sig(1)}>m_{\sig(2)}>m_{\sig(3)}>0}
		+
		\Sm{\tau\in\gpAlt{3}} \Sm{m_{\tau(1)}=m_{\tau(2)}>m_{\tau(3)}>0}
		\lnAHs{113}
		+
		\Sm{\tau\in\gpAlt{3}} \Sm{m_{\tau(1)}>m_{\tau(2)}=m_{\tau(3)}>0}
		+
		\Sm{m_1=m_2=m_3>0}
	}
	we can similarly see that	
	\envHLineCm
	{
		\lgP{k_1}{z_1}\lgP{k_2}{z_2}\lgP{k_3}{z_3}
	}
	{
		\Sm{\sig\in\gpSym{3}} \lgP{k_{\sig(1)},k_{\sig(2)},k_{\sig(3)}}{z_{\sig(1)},z_{\sig(1)}z_{\sig(2)},z_{\sig(1)}z_{\sig(2)}z_{\sig(3)}  }
		\lnAHs{0}
		+
		\Sm{\tau\in\gpAlt{3}} \lgP{k_{\tau(1)}+k_{\tau(2)},k_{\tau(3)}}{z_{\tau(1)}z_{\tau(2)},z_{\tau(1)}z_{\tau(2)}z_{\tau(3)}}
		\lnAHs{0}
		+
		\Sm{\tau\in\gpAlt{3}} \lgP{k_{\tau(1)},k_{\tau(2)}+k_{\tau(3)}}{z_{\tau(1)}, z_{\tau(1)}z_{\tau(2)}z_{\tau(3)}}
		\lnAHs{20}
		+
		\lgP{k_1+k_2+k_3}{z_{1}z_{2}z_{3}}
	}
	which with $z_1=z_2=z_3=z$ gives \refEq{2.1_Lem1Eq2}.
}
%Z
We now prove \refProp{2.1_Prop1}.
%A
\envProof[\refProp{2.1_Prop1}]{
We begin by verifying the first equation of \refEq{2.1_Prop1i_Eq1}.
For complex numbers  $X_1,X_2,X_3$,
	we find from the partial fraction expansion
	\envM
	{
		\opF[s]{1}{XY}
	}
	{
		\opF[s]{1}{(X+Y)Y} + \opF[s]{1}{(Y+X)X}
	}
	that
	\envOTLineCm[2.1_Prop1P_Eq0]
	{
		\opF{1}{X_{12}X_2X_3}
	}
	{
		\opF{1}{X_{312}X_{12}X_2} + \opF{1}{X_{123}X_2X_3}
	}
	{
		\Sm{\sig\in\lsetS{321}} \opF{1}{ X_{\sig(1)\sig(2)\sig(3)}  X_{\sig(2)\sig(3)} X_{\sig(3)} }
	}		
	where we make use of the definitions of \refEq{1_Def_ijk} for $X_{ij}$ and $X_{ijk}$.
We set the left- and the right-hand sides of \refEq{2.1_Prop1P_Eq0} equal to $\Fc{L}{X_1,X_2,X_3}$ and $\Fc{R}{X_1,X_2,X_3}$, respectively.
Since we see from 
	\envM{ \opF[s]{1}{(1-X)} }{ \Sm{l\geq1} X^{l-1} }
	that
	\envHLineCSCm[2.1_Prop1P_Eq1]
	{
		\Sm{m_1\geq1} \opF{z^{m_1}}{m_1-x_1t}
	}
	{
		\SmT{l=1}{\infty} x_1^{l-1} \lgP{l}{z} t^{l-1} 
	}
	{
		\Sm{m_1,m_2\geq1} \opF{z^{m_1+m_2}}{ (m_{12}-x_{12}t)(m_2-x_2t)}
	}
	{
		\SmT{l=2}{\infty} \gfcDPL{l}{x_{12},x_{2}}{z,z} t^{l-2}
	}
	{
		\Sm{m_1,m_2,m_3\geq1} \opF{z^{m_1+m_2+m_3}}{ (m_{123}-x_{123}t)(m_{23}-x_{23}t)(m_3-x_3t)}
	}
	{
		\SmT{l=3}{\infty} \gfcTPL{l}{x_{123},x_{23},x_{3}}{z,z,z} t^{l-3}
	}
	we have
	\envMLineCm
	{
		\Sm{m_1,m_2,m_3\geq1} z^{m_1+m_2+m_3}\Fc{L}{m_1-x_1t,m_2-x_2t,m_3-x_3t}
		\\
	}
	{
		\Sm{l\geq3} t^{l-3} \tpSm{l_1 \geq 2, l_2 \geq 1}{l_1+l_2=l} \gfcDPL{l_1}{x_{12},x_{2}}{z,z}  x_3^{l_2-1} \lgP{l_2}{z}
	}
	and
	\envMLinePd
	{
		\Sm{m_1,m_2,m_3\geq1} z^{m_1+m_2+m_3}\Fc{R}{m_1-x_1t,m_2-x_2t,m_3-x_3t}
		\\
	}
	{
		\Sm{l\geq3} t^{l-3} \Sm{\sig\in\lsetS{321}} \gfcTPL{l}{x_{\sig(1)\sig(2)\sig(3)},x_{\sig(2)\sig(3)},x_{\sig(3)}}{z,z,z}
	}
Calculating the coefficient of $t^{l-3}$ yields the first equation of \refEq{2.1_Prop1i_Eq1}.

Next, we prove the second equation of \refEq{2.1_Prop1i_Eq1}.
We see from \refEq{2.1_Lem1Eq1} that
	\envLineFPd
	{
		\tpSm{l_1 \geq 2, l_2 \geq 1}{l_1+l_2=l} \gfcDPL{l_1}{x_{1},x_{2}}{z,z^2}  x_3^{l_2-1} \lgP{l_2}{z}
	}
	{
		\tpSm{k_1,k_2,k_3 \geq 1}{k_1+k_2+k_3=l} x_1^{k_1-1}x_2^{k_2-1}x_3^{k_3-1} \lgP{k_1,k_2}{z,z^2} \lgP{k_3}{z}
	}
	{
		\tpSm{k_1,k_2,k_3 \geq 1}{k_1+k_2+k_3=l} x_1^{k_1-1}x_2^{k_2-1}x_3^{k_3-1}
		\lnAHs{50}
		\times
		\bkS[G]{ 
			\Sm{\sig\in\lsetS{321}} \lgP{k_{\sig(1)},k_{\sig(2)},k_{\sig(3)}}{z,z^2,z^3} 
			+ 
			\lgP{k_1+k_3,k_2}{z^2,z^3} 
			+ 
			\lgP{k_1,k_2+k_3}{z,z^3} 
		}
	}
	{
		\Sm{\sig\in\lsetS{321}} \gfcTPL{l}{x_{\sig(1)},x_{\sig(2)},x_{\sig(3)}}{z,z^2,z^3}
		\lnAHs{0}
		+
		\tpSm{k_1,k_2,k_3 \geq 1}{k_1+k_2+k_3=l} x_1^{k_1-1}x_2^{k_2-1}x_3^{k_3-1}
		\bkS[B]{ \lgP{k_1+k_3,k_2}{z^2,z^3} + \lgP{k_1,k_2+k_3}{z,z^3} }
	}
Since 
	\envMCm{
		\SmT{j=1}{m-1}X^{j-1}Y^{m-1-j}
	}{
		\opF[s]{(X^{m-1}-Y^{m-1})}{(X-Y)}
	}
	it follows that 
	\envLineFCm
	{
		\tpSm{k_1,k_2,k_3 \geq 1}{k_1+k_2+k_3=l} x_1^{k_1-1}x_2^{k_2-1}x_3^{k_3-1} \lgP{k_1+k_3,k_2}{z^2,z^3}
	}
	{
		\tpSm{l_1\geq2,l_2\geq1}{l_1+l_2=l} \bkR{ \SmT{j=1}{l_1-1}x_1^{j-1}x_3^{l_1-1-j} } x_2^{l_2-1} \lgP{l_1,l_2}{z^2,z^3}
	}
	{
		\tpSm{l_1, l_2\geq1}{l_1+l_2=l} \opF{x_1^{l_1-1}-x_3^{l_1-1}}{x_1-x_3} x_2^{l_2-1} \lgP{l_1,l_2}{z^2,z^3}
	}
	{
		\Sm{\sig\in\lsetS{31}} \opF{\gfcDPL{l}{x_{\sig(1)},x_{2}}{z^2,z^3}}{x_{\sig(1)} -  x_{\sig(3)}} 
	}
	and similarly that
	\envHLinePd
	{
		\tpSm{k_1,k_2,k_3 \geq 1}{k_1+k_2+k_3=l} x_1^{k_1-1}x_2^{k_2-1}x_3^{k_3-1} \lgP{k_1,k_2+k_3}{z,z^3}
	}
	{
		\Sm{\sig\in\lsetS{32}} \opF{\gfcDPL{l}{x_{1},x_{\sig(2)}}{z, z^3}}{x_{\sig(2)} -  x_{\sig(3)}} 	
	}
These yield the second equation of \refEq{2.1_Prop1i_Eq1}.

We can similarly prove the first equation of \refEq{2.1_Prop1ii_Eq1} by using the partial fraction expansion 
	\envOTLineCm
	{
		\opF{1}{X_1X_2X_3}
	}
	{
		\opF{1}{X_{12}X_2X_3} + \opF{1}{X_{21}X_1X_3}
	}
	{
		\Sm{\sig\in\gpSym{3}} \opF{1}{ X_{\sig(1)\sig(2)\sig(3)}  X_{\sig(2)\sig(3)} X_{\sig(3)} }
	}
	which follows from \refEq{2.1_Prop1P_Eq0}.
We can also prove the second equation of \refEq{2.1_Prop1ii_Eq1} by using \refEq{2.1_Lem1Eq2} and 
	\envHLinePd[2.1_Prop1PEq2]
	{
		\tpSm{k_1,k_2,k_3 \geq 1}{k_1+k_2+k_3=l} x_1^{k_1-1}x_2^{k_2-1}x_3^{k_3-1}
	}
	{
		\Sm{\tau\in\gpAlt{3}}  \bkR{\ip{x}{\tau(1)}{l-1} \nPdT{i=2}{3} \opF{1}{\pw{x_{\tau(1)}-x_{\tau(i)}}{1}}} 
	}
The equation \refEq{2.1_Prop1PEq2} is derived as follows:
We first verify that
	\envHLine
	{
		\Sm{\tau\in\gpAlt{3}} \opF{1}{X_{\tau(1)}} \PdT{i=2}{3} \opF{1}{X_{\tau(1)} -X_{\tau(i)} }
	}
	{
		\opF{1}{X_1X_2X_3}
	}
	by a direct calculation.
We next substitute $(x_1-t,x_2-t,x_3-t)$ for $(X_1,X_2,X_3)$ in this equation 
	and differentiate it with respect to $t$ $(l-3)$ times and evaluate it at $t=0$.
By replacing $x_i$ with $x_i^{-1}$ for $i=1,2,3$,
	we obtain \refEq{2.1_Prop1PEq2}.
}%Z
%A
%:2.1_Rem1
\begin{remark}\label{2.1_Rem1}
A coefficient of a $x_1^p x_2^q x_3^r$ term of the second equation of \refEq{2.1_Prop1i_Eq1} \lrpTx{\refEq{2.1_Prop1ii_Eq1}} 
	gives
	a harmonic relation \refEq{2.1_Lem1Eq1} \lrpTx{\refEq{2.1_Lem1Eq2}} and vice versa.
Thus, the second equations of \refEq{2.1_Prop1i_Eq1} and \refEq{2.1_Prop1ii_Eq1}  
	correspond to the harmonic relations for triple polylogarithms 
	involving two-fold and three-fold products, respectively.
	
We also see that
	the first equations of \refEq{2.1_Prop1i_Eq1} and \refEq{2.1_Prop1ii_Eq1} respectively 
	correspond to the shuffle relations for triple polylogarithms involving two-fold and three-fold products as follows.
	
Let $(r,q,p)$ be a three-tuple of positive integers.
By replacing $x_1$ with $x_1-x_2$ in the first equation of \refEq{2.1_Prop1i_Eq1}
	and calculating the coefficient of $ x_1^{r-1}x_2^{q-1}x_3^{p-1} $,
	we obtain
	\envHLineCm
	{
		\lgP{r,q}{z,z} \lgP{p}{z}
	}
	{
		\tpSm{l_1,l_2\geq1}{l_1+l_2=l-q} \binom{l_1-1}{p-1} \lgP{l_1,l_2,q}{z,z,z} 
		\lnAHs{0}
		+
		\tpSm{l_1,l_2,l_3\geq1}{l_1+l_2+l_3=l} \binom{l_1-1}{r-1} \bkS[G]{ \binom{l_2-1}{q-1} + \binom{l_2-1}{q-l_3} } \lgP{l_1,l_2,l_3}{z,z,z} 
	}
	where $\binom{m}{n} = 0$ if $m<n$ or $n<0$.
By replacing the multiple polylogarithms $\lgP{k_1,\ldots, k_n}{z,\ldots, z}$ with the words $z_{k_1}\cdots z_{k_n}$,
	this equation becomes the shuffle relation derived from (10) and (29) in \cite{KMT11},
	which gives the correspondence of the first equation of \refEq{2.1_Prop1i_Eq1} to the shuffle relations for triple polylogarithms involving two-fold product.
	
We next show the correspondence of the first equation of \refEq{2.1_Prop1ii_Eq1} 
	to the shuffle relations for triple polylogarithms involving three-fold product.
We can obtain
	\envHLine[2.1_Rem1_Eq1]
	{
		\tpSm{l_1,l_2\geq1}{l_1+l_2=l} x_1^{l_1-1}x_2^{l_2-1} \lgP{l_1}{z}\lgP{l_2}{z}
	}
	{
		\gfcDPL{l}{x_{12},x_{2}}{z,z} + \gfcDPL{l}{x_{21},x_{1}}{z,z}
	}
	using the second equation of \refEq{2.1_Prop1P_Eq1} and the partial fraction expansion 
	\envMPd{ \opF[s]{1}{XY} }{ \opF[s]{1}{(X+Y)Y} + \opF[s]{1}{(Y+X)X} }
We find that \refEq{2.1_Rem1_Eq1} corresponds to the shuffle relations for double polylogarithms;
	that is, 
	a coefficient of a $x_1^{p}x_2^{q}$ term in this equation gives a shuffle relation \cite[(24)]{KMT11} and vise versa.
Since 
	\envHLineCm
	{
		\lfcSPLi{l}{3}{x_1,x_2,x_3}{z}
	}
	{
		\tpSm{l_1\geq2, l_2 \geq 1}{l_1+l_2=l} 
		\bkR[B]{ \gfcDPL{l_1}{x_{12},x_{2}}{z,z} + \gfcDPL{l_1}{x_{21},x_{1}}{z,z} } x_3^{l_2-1} \lgP{l_2}{z}
	}
	the required correspondence is thus reduced to the case of \refEq{2.1_Prop1i_Eq1}.
\end{remark}%Z

\subsection{\sSectTwoT} \label{sSectTwoT} 
The asymptotic expansions of multiple polylogarithms have been given previously in \cite{IKZ06}.
In this  subsection,  
	using these expansions,
	we calculate constant terms of asymptotic expansions of the functions
	\envPLineCm[2.2_Plane_Fc1]
	{\lnA
		\gfcDPL{l}{x_1,x_2}{z,z}, 	\quad	\gfcDPL{l}{x_1,x_2}{z^m,z^n} \;\; (1\leq m< n\leq 3),
		\lnAH[]
		\gfcTPL{l}{x_{ 1 },x_{ 2 },x_{3}}{z,z,z}, 	\quad		\gfcTPL{l}{x_{ 1 },x_{ 2 },x_{3}}{z,z^2,z^3}
		\nonumber
	}
	which appear in \refProp{2.1_Prop1}.
	
In order to introduce the asymptotic expansions of multiple polylogarithms,
	we review the algebraic setup given by Hoffman \cite{Hoffman97}.
Let $\lalgFZV := \setQ \bkA{x,y}$ be the non-commutative polynomial algebra over the rational numbers
	in the two indeterminate letters $x$ and $y$,
	and let $\lalgFZV[1]$ and $\lalgFZV[0]$ be its subalgebras $\setQ+\lalgFZV y$ and $\setQ+x \lalgFZV y$, respectively.
We define the \emph{shuffle product} $\sh$ on $\lalgFZV$ inductively as follows:
	\envPLine
	{
		&1 \lsh w \lnP{=} w \lsh 1 \lnP{=} w,&
		\\
		&uw_1 \lsh vw_2 \lnP{=} u(w_1\lsh vw_2) + v(uw_1 \lsh w_2)&
	}
	for words $w,w_1,w_2\in\lalgFZV$ and $u,v\in\bkB{x,y}$ and extend it by $\setQ$-linearity. 
This product gives $\lalgFZV$ the structure of a commutative $\setQ$-algebra \cite{Reutenauer93},
	which we denote by $\lalgFZVs$.
The subspaces $\lalgFZV[1]$ and $\lalgFZV[0]$ also become subalgebras of $\lalgFZVs$ 
	and are denoted by $\lalgFZVs[1]$ and $\lalgFZVs[0]$, respectively.
Let  $\lmpZV[?]{} : \lalgFZV[0] \rightarrow \setR$ be the evaluation map defined in \cite[\S1]{IKZ06}; 
	that is,
	\envMDef
	{
		\lmpZV{w}
	}
	{
		\fcZ{l_1,\ldots,l_n}
	}
	for any word $w=x^{l_1-1}y\cdots x^{l_n-1}y\in\lalgFZV[0]$,
	and 
	let $\lmpZVs[?]{} : \lalgFZVs[1] \rightarrow \setR [T]$ be the algebra homomorphism defined in \cite[\S2]{IKZ06}.
We denote the image of the word $w=x^{l_1-1}y\cdots x^{l_n-1}y\in\lalgFZV[1]$ under the map $\lmpZVs[?]{}$ by
	\envHLineDefPd[2.2_Plain_DefEvMp]
	{
		\lmpZVi{l_1,\ldots,l_n}{T}
	}
	{
		\lmpZVs{w}
	}
This map describes the asymptotic properties of multiple polylogarithms (see \cite[p. 311 in \S2]{IKZ06}).
For any multi-index $(l_1,\ldots,l_n)$ of positive integers, 
	there exists a positive number $J>0$ such that
	\envHLineCm[2.2_Plain_AspPL]
	{\hspace{-20pt}
		\lgP{l_1,\ldots,l_n}{z,\ldots,z}
	}
	{
		\lmpZVi[b]{l_1,\ldots,l_n}{-\lgg[n]{1-z}} + \sbLandau[B]{(1-z) \pwR[b]{\lgg[n]{1-z}}{J}}	\quad	(z\nearrow1)
	}
	where $\sbLandau[]{}$ denotes the Landau symbol.	
(We see that \envM{ \lgP{l_1,\ldots,l_n}{z,\ldots,z} }{ \lgP{l_1,\ldots,l_n}{z} } by comparing \refEq{1_Def_lgP} and \cite[(2.4)]{IKZ06}.)

For any function $\Fc{f}{z}$ which has a polynomial $\Fc{P}{T}$ and a positive number $J>0$
	and 
	satisfies an asymptotic property of the form \refEq{2.2_Plain_AspPL},
	we denote the constant term of $\Fc{P}{T}$ or $\Fc{P}{0}$ by $\lfcCst{\Fc{f}{z}}$.
For example, 
	\envHLinePd[2.2_Plain_EqCtTmAndEvMp]
	{ 
		\lfcCst[b]{\lgP{l_1,\ldots,l_n}{z,\ldots,z}} 
	}
	{ 
		\lmpZVi{l_1,\ldots,l_n}{0} 
	}	
By \refEq{2.2_Plain_DefEvMp}, 
	the image of $\lgP{l_1,\ldots,l_n}{z,\ldots,z}$ under $\lfcCst[]{}$ can also be expressed as the composition of the evaluation map $\lmpZV[?]{}$ 
	and 
	the regularization map $\lmpRegT[?]{}:\lalgFZVs[1] \rightarrow \lalgFZVs[0] [T]$ with $T=0$ which is defined in \cite[\S3]{IKZ06},
	\envHLineCm[2.2_Plain_EqCtTmAndRegMp]
	{ 
		\lfcCst{\lgP{l_1,\ldots,l_n}{z,\ldots,z}} 
	}
	{ 
		\lmpZV{ \lmpReg{x^{l_1-1}y\cdots x^{l_n-1}y} }
	}
	where $\lmpReg[?]{} = \lmpRegT[?]{} \vert_{T=0}$.
	
For positive integers $m, n, n_1, n_2$ such that $m\geq3, n\geq4, n_1\geq2, n_2\geq1$, 
	we define the real values $\lmfcZ{1,m-1}, \lmfcZ{1,1,n-2}$, and $\lmfcZ{1,n_1,n_2}$ as
	\envHLineCSDefPd[2.2_Def_CtDTZV]
	{
		\lmfcZ{1,m-1}
	}
	{
		- \bkR[G]{ \tpSm{j_1\geq2,j_2 \geq1}{j_1+j_2=m} \fcZ{j_1,j_2} + \fcZ{m-1,1} }
	}
	{
		\lmfcZ{1,1,n-2}
	}
	{
		\tpSm{j_1\geq2,j_2,j_3\geq1}{j_1+j_2+j_3=n} \fcZ{j_1,j_2,j_3} + \tpSm{j_1\geq2,j_2 \geq1}{j_1+j_2=n-1} \fcZ{j_1,j_2,1}+ \fcZ{n-2,1,1}
	}
	{
		\lmfcZ{1,n_1,n_2}
	}
	{
		- \bkR[G]{ \tpSm{j_1\geq2,j_2\geq1}{j_1+j_2=n_1+1} \fcZ{j_1,j_2,n_2} + \tpSm{j_2,j_3\geq1}{j_2+j_3=n_2+1} \fcZ{n_1,j_2,j_3} + \fcZ{n_1,n_2,1} }
	}
We also define
	 $\lmfcZ{1,1}=\lmfcZ{1,1,1}=0$.
The values defined in \refEq{2.2_Def_CtDTZV} are equal to 
	$\lmpZVi{1,m-1}{0} = \lmpZV{ \lmpReg{y x^{m-2}y} } $, 
	$\lmpZVi{1,1,n-2}{0} = \lmpZV{ \lmpReg{y^2 x^{n-3}y} }$, 
	and $\lmpZVi{1,n_1,n_2}{0} = \lmpZV{ \lmpReg{y x^{n_1-1}y x^{n_2-1}y } }$, respectively,
	by virtue of \refEq{2.2_Plain_EqCtTmAndEvMp}, \refEq{2.2_Plain_EqCtTmAndRegMp}, and \refLem{2.2_Lem1} below;
	that is,
	the values defined in \refEq{2.2_Def_CtDTZV} are the regularized double and triple zeta values related to the shuffle regularization.
	
Hoffman's relations \cite[Theorem 5.1]{Hoffman92}
	with $(i_1,i_2)=(l_1,l_2)$ and $k=2$ are
	\envHLineCm[2.2_Plane_EqHoffman]
	{
		\fcZ{l_1+1,l_2} + \fcZ{l_1,l_2+1}
	}
	{
		\tpSm{j_1\geq2,j_2\geq1}{j_1+j_2=l_1+1} \fcZ{j_1,j_2,l_2} + \tpSm{j_2\geq2,j_3\geq1}{j_2+j_3=l_2+1} \fcZ{l_1,j_2,j_3}
	}
	where an empty sum is defined as $0$.
From \refEq{2.2_Plane_EqHoffman} and the sum formulas \refEq{1_Plane_EqSFofDZV} and \refEq{1_Plane_EqSFofTZV},
	we find the following simple expressions for those of \refEq{2.2_Def_CtDTZV},
	\envHLineCSPd[2.2_Plane_EqCtDTZV]
	{
		\lmfcZ{1,m-1}
	}
	{
		-\bkR[b]{\fcZ{m-1,1} + \fcZ{m}}
	}
	{
		\lmfcZ{1,1,n-2}
	}
	{
		\fcZ{n-2,1,1} + \fcZ{n-1,1} + \fcZ{n-2,2} + \fcZ{n}
	}
	{
		\lmfcZ{1,n_1,n_2}
	}
	{
		-\bkR[b]{ \fcZ{n_1,n_2,1} + \fcZ{n_1,1,n_2} + \fcZ{n_1+1,n_2} + \fcZ{n_1,n_2+1} }	
	}
	
By the use of the algebraic formula \cite[Proposition 8]{IKZ06} for the regularization map $\lmpRegT[?]{}$
	and  \refEq{2.2_Plain_EqCtTmAndRegMp}, 
	we can explicitly calculate $\lfcCst{\Fc{f}{z}}$
	for some functions $f(z)$ expressed in terms of double and triple polylogarithms.
%A
%:2.2_Lem1
\begin{lemma}\label{2.2_Lem1}
Let $k,k_1,k_2$ be positive integers, and $\delta_{m,n}$ be the Kronecker delta function.
\\{\bf (i)} 
We have
	\envHLineCFPd[2.2_Lem1i_Eq1]
	{
		\lfcCst{\lgP{1,1}{z,z} }
	}
	{
		0
	}
	{
		\lfcCst{\lgP{1,1}{z,z^2} }
	}
	{
		- \opF{1}{2} \fcZ{2}
	}
{\bf (ii)} 
If $k \geq 3$, then for positive integers $m,n$ such that $1\leq m<n\leq 3$, we have 
	\envHLineCFPd[2.2_Lem1ii_Eq1]
	{
		\lfcCst{\lgP{1,k-1}{z,z} }
	}
	{
		\lmfcZ{1,k-1}
	}
	{
		\lfcCst{ \lgP{1,k-1}{z^m,z^n} }
	}
	{
		\lmfcZ{1,k-1} - \delta_{m,2} \fcZ{k-1}
	}
{\bf (iii)}
If $k\geq4$, then we have
	\envHLineCFPd[2.2_Lem1iii_Eq1]
	{
		\lfcCst{ \lgP{1,1,k-2}{z,z,z}  }
	}
	{
		\lmfcZ{1,1,k-2}
	}
	{
		\lfcCst{ \lgP{1,1,k-2}{z,z^2,z^3}  }
	}
	{
		\lmfcZ{1,1,k-2} - \opF{1}{2} \fcZ{k-2}\fcZ{2}
	}
{\bf (iv)}
If $k_1\geq2$, then we have
	\envOTLinePd[2.2_Lem1iv_Eq1]
	{
		\lfcCst{ \lgP{1,k_1,k_2}{z,z,z}  }
	}
	{
		\lfcCst{ \lgP{1,k_1,k_2}{z,z^2,z^3}  }
	}
	{
		\lmfcZ{1,k_1,k_2}
	}
\end{lemma}
%Z
%A
\envProof{
Using the partial fraction expansion 
	\envM{ \opF[s]{1}{(m_1+m_2)m_2}+\opF[s]{1}{(m_2+m_1)m_1} }{ \opF[s]{1}{m_1m_2} }
	and the summation decomposition 
	\envM{ \Sm{m_1,m_2>0} }{ \Sm{m_1>m_2>0} + \Sm{m_2>m_1>0} + \Sm{m_1=m_2>0} },
	we obtain
	\envOTLineCm[2.2_Lem1P_Eq1]
	{
		2\lgP{1,1}{z,z}
	}
	{
		 \lgP{1}{z}^2
	}
	{
		2 \lgP{1,1}{z,z^2} + \lgP{2}{z^2}
	}
	which with \envM{ \lgP{1}{z} }{ - \lgg[a]{1-z} } yields the equations of \refEq{2.2_Lem1i_Eq1}.

The following appears in the proof of \cite[\refPropF{8}]{IKZ06}:
	\envHLineCm[2.2_Lem1P_Eq2]
	{
		\lmpRegT{y^m w_0}
	}
	{
		\SmT{j=0}{m} \mo^{j} x (y^j \lsh w_0') \opF{T^{m-j}}{(m-j)!}
	}
	where $w_0\in\lalgFZV[0]$ and $w_0'\in\lalgFZV[1]$ such that $w_0=x w_0'$.
Since $\lmpReg[?]{} = \lmpRegT[?]{} |_{T=0}$,
	we obtain
	\envM{ \lmpZV{ \lmpReg{y^m w_0 } }  }{ \lmpZV{ \mo^m x (y^m \lsh w_0') } },
	which with \refEq{2.2_Plain_EqCtTmAndRegMp} yields 
	\envHLineCSPd
	{
		\lfcCst{\lgP{1,k-1}{z,z}} 
	}
	{
		- \lmpZV{x (y \lsh x^{k-3}y) }
	}
	{
		\lfcCst{\lgP{1,1,k-2}{z,z,z}} 
	}
	{
		\lmpZV{x (y^2 \lsh x^{k-4}y) }
	}
	{
		\lfcCst{\lgP{1,k_1,k_2}{z,z,z}} 
	}
	{
		 -\lmpZV{x (y \lsh x^{k_1-2}y x^{k_2-1}y ) }
	}
By \refEq{2.2_Def_CtDTZV} and 
	\envHLineCSCm
	{
		x (y \lsh x^{k-3}y)
	}
	{
		\tpSm{j_1\geq2,j_2 \geq1}{j_1+j_2=k}x^{j_1-1} y x^{j_2-1}y + x^{k-2}y^2
	}
	{
		x (y^2 \lsh x^{k-4}y)
	}
	{
		\tpSm{j_1\geq2,j_2,j_3\geq1}{j_1+j_2+j_3=k} x^{j_1-1} y x^{j_2-1}y x^{j_3-1} y 
		\lnAHs{0}
		+ 
		\tpSm{j_1\geq2,j_2 \geq1}{j_1+j_2=k-1}x^{j_1-1} y x^{j_2-1}y^2 + x^{k-3}y^3
	}
	{
		x (y \lsh x^{k_1-2}y x^{k_2-1}y)
	}
	{
		\tpSm{j_1\geq2,j_2\geq1}{j_1+j_2=k_1+1} x^{j_1-1} y x^{j_2-1}y x^{k_2-1} y
		\lnAHs{0}
		+
		\tpSm{j_2,j_3\geq1}{j_2+j_3=k_2+1} x^{k_1-1} y x^{j_2-1}y x^{j_3-1} y
		+
		x^{k_1-1} y x^{k_2-1}y^2
	}
	we verify the first equations of \refEq{2.2_Lem1ii_Eq1} and \refEq{2.2_Lem1iii_Eq1}, 
	and $\lfcCst{ \lgP{1,k_1,k_2}{z,z,z}  }=\lmfcZ{1,k_1,k_2}$.
	
We will next prove the second equation of \refEq{2.2_Lem1ii_Eq1}.
If $z$ is a real number such that $0<z<1$,
	we see by the assumption $k\geq3$ that
	\envHLineFiCm[2.2_Lem1P_Eq3]
	{
		\lgP{1,k-1}{z,z} - \lgP{1,k-1}{z,z^{3}}
	}
	{
		\Sm{m_1>m_2>0} \opF{z^{m_1} - z^{m_1+2m_2}}{m_1m_2^{k-1}}
	}
	{
		\Sm{m_1>m_2>0} \opF{z^{m_1}(1-z^{2 m_2})}{m_1m_2^{k-1}}
	}
	{
		(1-z) \Sm{m_1>m_2>0} \opF{z^{m_1}(\sqSm{1+z}{z^{2m_2-1})}} {m_1m_2^{k-1}}
	\lnAHP{\leq}
		2 (1-z) \Sm{m_1>m_2>0} \opF{z^{m_1}}{m_1\pw{m_2}{k-2}}		
	\lnAHP{\leq}
		2 (1-z) \Sm{m_1>m_2>0} \opF{z^{m_1}} {m_1 m_2}		
	}
	{
		2 (1-z)  \lgP{1,1}{z,z}
	}
	which with the first equation of \refEq{2.2_Lem1ii_Eq1}
	proves the second one for $(m,n)=(1,3)$, because of \refEq{2.2_Lem1P_Eq1} and that \envM{ \lgP{1}{z} }{ - \lgg[a]{1-z} }.
The case of $(m,n)=(1,2)$ can be proved from those of $(m,n)=(1,1)$ and $(1,3)$
	and the fact that $ \lgP{1,k-1}{z,z^3}  <  \lgP{1,k-1}{z,z^2} < \lgP{1,k-1}{z,z}$ for any real number $z$ such that $0<z<1$.
For the remaining case, $(m,n)=(2,3)$,
	we see that
	\envHLineFPd
	{
		\lgP{1,k-1}{z,z^3} - \lgP{1,k-1}{z^2,z^{3}}
	}
	{
		\Sm{m_1>m_2>0} \opF{z^{m_1+2m_2} - z^{2m_1+m_2}}{m_1m_2^{k-1}}
	}
	{
		\Sm{m_1>m_2>0} \opF{z^{m_1+2m_2}(1-z^{m_1- m_2})}{m_1m_2^{k-1}}
	}
	{
		(1-z) \Sm{m_1>m_2>0} \opF{z^{m_1+2m_2}(\sqSm{1+z}{z^{m_1-m_2-1})}} {m_1m_2^{k-1}}
	\lnAHP{\leq}
		(1-z) \Sm{m_1>m_2>0} \opF{z^{m_1+m_2}(m_1-m_2)}{m_1m_2^{k-1}}		
	}
We allow $l,l_1$, or $l_2$ to be equal to $0$ in the definitions of the polylogarithms $\lgP{l}{z}$ and $\lgP{l_1,l_2}{z_1,z_2}$,
	which are well-defined because $\abs{z},\abs{z_i}<1$.
Note that \envM{\lgP{0}{z}}{z/(1-z)}.
Since
	\envHLineFiCm
	{
		\lgP{k-1,0}{z,z^2}
	}
	{
		\Sm{m_1>m_2>0} \opF{z^{m_1+m_2}}{ m_1^{k-1} }
	}
	{
		\Sm{m_1>0} \opF{z^{m_1}}{m_1^{k-1}} \opF{1-z^{m_1}}{1-z} - \lgP{k-1}{z}
	}
	{
		\opF{1}{1-z} \bkR{ \lgP{k-1}{z}-\lgP{k-1}{z^2} } - \lgP{k-1}{z}
	}
	{
		\lgP{0}{z} \lgP{k-1}{z} - \opF{1}{z} \lgP{0}{z}\lgP{k-1}{z^2}
	}
	we have
	\envLineFPd
	{
		\Sm{m_1>m_2>0} \opF{ z^{m_1+m_2}(m_1-m_2) }{ m_1m_2^{k-1} }		
	}
	{
		\bkR[G]{ \Sm{m_1,m_2>0} - \Sm{m_2>m_1>0} } \bkR[G]{ \opF{ z^{m_1+m_2} }{ m_2^{k-1} }-\opF{ z^{m_1+m_2} }{ m_1\pw{m_2}{k-2} } }
	}
	{
		\lgP{0}{z}\lgP{k-1}{z} - \lgP{1}{z}\lgP{k-2}{z} - \lgP{k-1,0}{z,z^2} + \lgP{k-2,1}{z,z^2}
	}
	{
		- \lgP{1}{z}\lgP{k-2}{z} + \opF{1}{z} \lgP{0}{z}\lgP{k-1}{z^2} + \lgP{k-2,1}{z,z^2}
	}
Therefore,
	\envHLineCmPt{\leq}
	{
		\lgP{1,k-1}{z,z^3} - \lgP{1,k-1}{z^2,z^{3}} - \lgP{k-1}{z^2}
	}
	{
		(1-z) \bkR[b]{ - \lgP{1}{z}\lgP{k-2}{z} + \lgP{k-2,1}{z,z^2} }
	}
	which gives the second equation of \refEq{2.2_Lem1ii_Eq1} for $(m,n)=(2,3)$,
	and so
	we have completed the proof of \refEq{2.2_Lem1ii_Eq1}.
	
We will next prove \refEq{2.2_Lem1iv_Eq1}.
We see from the first line of \refEq{2.1_Prop1ii_Eq1} with $l=3$ that
	\envOTLinePd
	{ 
		6\lgP{1,1,1}{z,z,z} 
	}
	{ 
		\lgP{1}{z}^3 
	}
	{ 
		\pwR{- \lgg[a]{1-z}}{3} 
	}
By the assumption $k_1\geq2$,
	a calculation similar to \refEq{2.2_Lem1P_Eq3} with this equation and \refEq{2.2_Lem1P_Eq1} gives 
	\envHLineFPdPte
	{
		\lgP{1,k_1,k_2}{z,z,z} - \lgP{1,k_1,k_2}{z,z^2,z^3}
	}{=}
	{
		\Sm{m_1>m_2>m_3>0} \opF{z^{m_1}(1-z^{m_2+m_3})}{m_1m_2^{k_2} m_3^{k_3}}
	}{\leq}
	{
		(1-z) \bkR{ \lgP{1,1,1}{z,z,z} + \lgP{1,1}{z,z} }	
	}{=}
	{
		(1-z) \bkR[g]{ \opF{\pwR{- \lgg[a]{1-z}}{3} }{6} + \opF{\pwR{- \lgg[a]{1-z}}{2} }{2} }
	}
Thus, 
	we see that 
	\envM{
		\lfcCst{ \lgP{1,k_1,k_2}{z,z,z} }
	}{
		\lfcCst{ \lgP{1,k_1,k_2}{z,z^2,z^3} }
	}
	and so 
	\refEq{2.2_Lem1iv_Eq1} follows
	since we have already proved that $\lfcCst{ \lgP{1,k_1,k_2}{z,z,z}  }=\lmfcZ{1,k_1,k_2}$.

For a  proof of the second equation of \refEq{2.2_Lem1iii_Eq1}, which is the final task, 
	we need the following identity derived from \refEq{2.1_Lem1Eq1} with $(k_1,k_2,k_3)=(1,k-2,1)$:
	\envLinePd
	{
		\lgP{1,k-2}{z,z^2}\lgP{1}{z}
	}
	{
		2\lgP{1,1,k-2}{z,z^2,z^3} + \lgP{1,k-2,1}{z,z^2,z^3} + \lgP{1,k-1}{z,z^3} + \lgP{2,k-2}{z^2,z^3}
	}
By the condition $k\geq4$, 
	this identity, along with \refEq{2.2_Plane_EqCtDTZV}, \refEq{2.2_Lem1ii_Eq1}, and \refEq{2.2_Lem1iv_Eq1}, yields 
	\envLineFCm
	{
		2\lfcCst[b]{ \lgP{1,1,k-2}{z,z^2,z^3} }
	}
	{
		-
		\lfcCst[b]{ \lgP{1,k-2,1}{z,z^2,z^3} }
		- 
		\lfcCst[b]{ \lgP{1,k-1}{z,z^3} }
		-
		\lfcCst[b]{ \lgP{2,k-2}{z^2,z^3} }
	}
	{
		\bkR[b]{ 2\fcZ{k-2,1,1} + \fcZ{k-1,1} + \fcZ{k-2,2} } + \bkR[b]{\fcZ{k-1,1} + \fcZ{k}} - \fcZ{2,k-2}
	}
	{
		2 \lmfcZ{1,1,k-2} -\fcZ{k-2}\fcZ{2}
	}
	which proves the second equation of \refEq{2.2_Lem1iii_Eq1}.
}
%Z
By \refLem{2.2_Lem1}, 
	we can calculate the constant terms of asymptotic expansions of the functions \refEq{2.2_Plane_Fc1} (see \refProp{2.2_Prop1} below).
We omit the proof of the proposition since it is obvious.
To describe the constant terms concisely, 
	we prepare the generating functions 
	$\lmgfcDZV{l}{x_1,x_2}$ and $\lmgfcTZV{l}{x_1,x_2,x_3}$
	of the regularized double and triple zeta values.
By \refEq{2.2_Plain_EqCtTmAndEvMp} and \refLem{2.2_Lem1},
	they are defined by
	\envHLineCFDefPd[2.2_Plane_DefMGFC]
	{
		\lmgfcDZV{l}{x_1,x_2}
	}
	{
		\gfcDZV{l}{x_1,x_2} + x_2^{l-2} \lmfcZ{1,l-1}
	}
	{
		\lmgfcTZV{l}{x_1,x_2,x_3}
	}
	{
		\gfcTZV{l}{x_1,x_2,x_3} + \tpSm{l_2,l_3\geq1}{l_2+l_3=l-1} x_2^{l_2-1}x_3^{l_3-1} \lmfcZ{1,l_2,l_3} 
	}
%A
%:2.2_Prop1
\begin{proposition}\label{2.2_Prop1}
Let $l$ be a positive integer.
\\{\bf (i)} 
We have
	\envHLineCFPd[2.2_Prop1i_Eq1]
	{
		\lfcCst{ \gfcDPL{2}{x_1,x_2}{z,z} }
	}
	{
		0
	}
	{
		\lfcCst{ \gfcDPL{2}{x_1,x_2}{z,z^2} }
	}
	{
		- \opF{1}{2} \fcZ{2}
	}
{\bf (ii)} 
If $l\geq3$, then for positive integers $m,n$ such that $1\leq m<n\leq 3$, we have 
	\envHLineCFPd[2.2_Prop1ii_Eq1]
	{
		\lfcCst{ \gfcDPL{l}{x_1,x_2}{z,z} }
	}
	{
		\lmgfcDZV{l}{x_1,x_2}
	}
	{
		\lfcCst{ \gfcDPL{l}{x_1,x_2}{z^m,z^n} }
	}
	{
		\lmgfcDZV{l}{x_1,x_2} - \delta_{m,2}x_2^{l-2}\fcZ{l-1}
	}
{\bf (iii)} 
If $l\geq4$, then we have
	\envHLineCFPd[2.2_Prop1iii_Eq1]
	{
		\lfcCst{ \gfcTPL{l}{x_{ 1 },x_{ 2 },x_{3}}{z,z,z} }
	}
	{
		\lmgfcTZV{l}{x_1,x_2,x_3}
	}
	{
		\lfcCst{ \gfcTPL{l}{x_{ 1 },x_{ 2 },x_{3}}{z,z^2,z^3} }
	}
	{
		\lmgfcTZV{l}{x_1,x_2,x_3} - \opF{x_3^{l-3}}{2}\fcZ{l-2}\fcZ{2}
	}
\end{proposition}
%Z	

\section{\sectThree} \label{sectThree}
In this section, 
	we prove \refThm{1_Thm1}. 
More precisely, 
	we prove 
	the two formulas \refEq{1_Thm1i_Eq1} and \refEq{1_Thm1ii_Eq1} for the generating functions $\gfcDZV{l}{x_1,x_2}$  and $\gfcTZV{l}{x_1,x_2,x_3}$ 
	by the use of the two classes $\lclsEDS{l}{2,1}$ and $\lclsEDS{l}{1,1,1}$ of extended double shuffle relations for triple zeta values.

The class $\lclsEDS{l}{2,1}$ \lrpTx{$\lclsEDS{l}{1,1,1}$ }
	is expressed as \refEq{3_Prop1i_Eq1} \lrpTx{\refEq{3_Prop1ii_Eq1}} below 
	in the sense that 
	the coefficients of $x_1^{p-1}x_2^{q-1}x_3^{r-1}$ of \refEq{3_Prop1i_Eq1} \lrpTx{\refEq{3_Prop1ii_Eq1}}
	are equal to the extended double shuffle relations of $\lclsEDS{l}{2,1}$ \lrpTx{$\lclsEDS{l}{1,1,1}$}.
It should be noted that
	the coefficients except that for $p=r=1$ (which is that of $x_2^{l-3}$) of \refEq{3_Prop1i_Eq1}
	are the original extended double shuffle relations $\loclsEDS{l}{3}$ for triple zeta values in \cite[Theorem 2(iv)]{IKZ06},
	\envHLineCm
	{
		\lmpZV{ \lmpReg{ w_1 \lsh w_0 - w_1 * w_0 } }
	}
	{
		0
	}
	where $w_1\in\lalgFZV[1]$ and $w_0\in\lalgFZV[0]$ 
	such that the pair $(d_1,d_0)$ of the depths of $w_1$ and $w_0$ is equal to $(2,1)$ or $(1,2)$.
On the other hand, 
	the coefficient of $x_2^{l-3}$ of \refEq{3_Prop1i_Eq1} is not stated in \cite[Theorem 2(iv)]{IKZ06} 
	and includes the extra value $-\fcZ{l-2}\fcZ{2}$
	which does not appear in any extended double shuffle relation of $\loclsEDS{l}{3}$ apparently.

%A
%:3_Prop1
\begin{proposition}\label{3_Prop1}
Let $l$ be an integer such that $l\geq4$, and $x_1,x_2,x_3$ be variables.
We define
	\envHLineCFDefPd
	{
		\lfcSPZ{l}{2}{x_1,x_2,x_3}
	}
	{
		\tpSm{l_1 \geq 3, l_2 \geq 2}{l_1+l_2=l} \lmgfcDZV{l_1}{x_{1},x_{2}}  x_3^{l_2-1} \fcZ{l_2}
	}
	{
		\lfcSPZ{l}{3}{x_1,x_2,x_3}
	}
	{
		\tpSm{l_1,l_2,l_3\geq2}{l_1+l_2+l_3=l} x_1^{l_1-1}x_2^{l_2-1}x_3^{l_3-1} \fcZ{l_1}\fcZ{l_2}\fcZ{l_3}
	}
\\{\bf (i)} 
We have
	\envHLineThPd[3_Prop1i_Eq1]
	{
		\lfcSPZ{l}{2}{x_1,x_2,x_3}
	}
	{
		\lmgfcTZV{l}{x_{13},x_{23},x_{3}} + \lmgfcTZV{l}{x_{13},x_{32},x_{2}} + \lmgfcTZV{l}{x_{31},x_{1},x_{2}}
	}
	{
		\Sm{\sig\in\lsetS{321}} \lmgfcTZV{l}{x_{\sig(1)},x_{\sig(2)},x_{\sig(3)}}
		\lnAHs{0}
		+
		\Sm{\sig\in\lsetS{31}} \opF{\lmgfcDZV{l}{x_{\sig(1)},x_{2}}}{x_{\sig(1)} -  x_{\sig(3)}} 
		+
		\Sm{\sig\in\lsetS{32}} \opF{\lmgfcDZV{l}{x_{1},x_{\sig(2)}}}{x_{\sig(2)} - x_{\sig(3)}} 
		-
		x_2^{l-3}\fcZ{l-2}\fcZ{2}
	}
{\bf (ii)} 
We have
	\envHLineThPd[3_Prop1ii_Eq1]
	{
		\lfcSPZ{l}{3}{x_1,x_2,x_3}
	}
	{
		\Sm{\sig\in\gpSym{3}} \lmgfcTZV{l}{x_{\sig(1)\sig(2)\sig(3)},x_{\sig(2)\sig(3)},x_{\sig(3)}}
	}
	{
		\Sm{\sig\in\gpSym{3}}
		\bkS[G]{
			\lmgfcTZV{l}{x_{\sig(1)},x_{\sig(2)},x_{\sig(3)}}
			+
			\opF{\lmgfcDZV{l}{x_{\sig(1)},x_{\sig(2)}}}{x_{\sig(1)} - \pw{x_{\sig(3)}}{1}} 
			+
			\opF{\lmgfcDZV{l}{x_{\sig(1)},x_{\sig(2)}}}{x_{\sig(2)} - \pw{x_{\sig(3)}}{1}} 
		}
		\lnAHs{0}
		+
		\Sm{\tau\in\gpAlt{3}}
		\bkS[G]{
			\bkR[g]{x_{\tau(1)}^{l-1} \PdT{i=2}{3} \opF{1}{x_{\tau(1)}-x_{\tau(i)}}} \fcZ{l}
			-
			x_{\tau(1)}^{l-3} \fcZ{l-2}\fcZ{2}
		}
	}
\end{proposition}
%Z
%A
\envProof{
We see from the first equation of \refEq{2.1_Prop1i_Eq1} and \refProp{2.2_Prop1} that
	\envHLineCm
	{
		\tpSm{l_1 \geq 3, l_2 \geq 2}{l_1+l_2=l} \lmgfcDZV{l_1}{x_{12},x_{2}}  x_3^{l_2-1} \fcZ{l_2}
	}
	{
		\Sm{\sig\in\lsetS{321}} \lmgfcTZV{l}{x_{ \sig(1)\sig(2)\sig(3) },x_{ \sig(2)\sig(3) },x_{\sig(3)}}
	}
	which verifies the first line of \refEq{3_Prop1i_Eq1} by replacing $x_1$ with $x_1-x_2$.
We also see from the second equation of \refEq{2.1_Prop1i_Eq1} and \refProp{2.2_Prop1} that
	\envLinePd
	{
		\tpSm{l_1 \geq 3, l_2 \geq 2}{l_1+l_2=l} \lmgfcDZV{l_1}{x_{1},x_{2}}  x_3^{l_2-1} \fcZ{l_2}
		-
		\opF{ x_3^{l-3} }{2}\fcZ{l-2}\fcZ{2}
	}
	{
		\Sm{\sig\in\lsetS{321}} \lmgfcTZV{l}{x_{\sig(1)},x_{\sig(2)},x_{\sig(3)}}
		-
		\opF{ 2x_2^{l-3}+x_3^{l-3} }{2}\fcZ{l-2}\fcZ{2}
		\lnAHs{0}
		+
		\Sm{\sig\in\lsetS{31}} \opF{\lmgfcDZV{l}{x_{\sig(1)},x_{2}} - x_2^{l-2} \fcZ{l-1}}{x_{\sig(1)} -  x_{\sig(3)}}
		+
		\Sm{\sig\in\lsetS{32}} \opF{\lmgfcDZV{l}{x_{1},x_{\sig(2)}}}{x_{\sig(2)} - x_{\sig(3)}} 
	}
Since 
	\envMCm[d]{ \Sm{\sig\in\lsetS{31}} \opF{x_2^{l-2} }{x_{\sig(1)} -  x_{\sig(3)}} \fcZ{l-1} }{ 0 }
	this yields the second line of \refEq{3_Prop1i_Eq1}.
We can similarly prove the equalities of \refEq{3_Prop1ii_Eq1} 
	by using \refEq{2.1_Prop1ii_Eq1} instead of \refEq{2.1_Prop1i_Eq1},
	and so omit their proofs here.
}
%Z

We prove \refThm{1_Thm1}.
%A
\envProof[\refThm{1_Thm1}]{
We see from \refEq{2.2_Def_CtDTZV} and \refEq{2.2_Plane_DefMGFC} that
	$\lmgfcTZV{l}{x_1,x_2,x_3}$ is expressed as
	\envLinePd
	{
		\lmgfcTZV{l}{x_1,x_2,x_3}
	}
	{
		\gfcTZV{l}{x_1,x_2,x_3} - \gfcTZV{l}{x_2,x_2,x_3} - \gfcTZV{l}{x_2,x_3,x_3} - \gfcTZV{l}{x_2,x_3,0} +x_3^{l-3} \lmfcZ{1,1,l-2}
	}
Thus, by \refEq{3_Prop1i_Eq1}, we have
	\envLineCm[3_Thm1P_Eq1]
	{
		\gfcTZV{l}{x_{13},x_{23},x_{3}} + \gfcTZV{l}{x_{13},x_{32},x_{2}} + \gfcTZV{l}{x_{31},x_{1},x_{2}}
		\lnAHs[]{0}
		-
		\Sm{\sig\in\lsetS{32} } \bkR{ 
			\gfcTZV{l}{x_{\sig(2)\sig(3)},x_{\sig(2)\sig(3)},x_{\sig(3)}} 
			+ 
			\gfcTZV{l}{x_{\sig(2)\sig(3)},x_{\sig(3)},x_{\sig(3)}} 
			+ 
			\gfcTZV{l}{x_{\sig(2)\sig(3)},x_{\sig(3)},0} 
		}
		\nonumber
	}
	{
		\gfcTZV{l}{x_{1},x_{2},x_{3}} + \gfcTZV{l}{x_{1},x_{3},x_{2}} + \gfcTZV{l}{x_{3},x_{1},x_{2}}
		\lnAHs{0}
		-
		\Sm{\sig\in\lsetS{32} } \bkR{ 
			\gfcTZV{l}{x_{\sig(2)},x_{\sig(2)},x_{\sig(3)}} 
			+ 
			\gfcTZV{l}{x_{\sig(2)},x_{\sig(3)},x_{\sig(3)}} 
			+ 
			\gfcTZV{l}{x_{\sig(2)},x_{\sig(3)},0} 
		}
		\lnAHs{0}
		+
		\Sm{\sig\in\lsetS{31}} \opF{\gfcDZV{l}{x_{\sig(1)},x_{2}}}{x_{\sig(1)} -  x_{\sig(3)}} 
		+
		\Sm{\sig\in\lsetS{32}} \opF{\gfcDZV{l}{x_{1},x_{\sig(2)}}}{x_{\sig(2)} - x_{\sig(3)}} 
%		%%
%		\lnAHs{30}
		+
		\opF{x_2^{l-2}-x_3^{l-2}}{x_2-x_3} \lmfcZ{1,l-1}
		-
		x_2^{l-3}\fcZ{l-2}\fcZ{2}
	}
	which with $x_1=0$ yields
	\envLinePd[3_Thm1P_Eq2]
	{
		\gfcTZV{l}{x_{3},x_{23},x_{3}} + \gfcTZV{l}{x_{3},x_{32},x_{2}}
		\lnAHs[]{0}
		-
		\Sm{\sig\in\lsetS{32} } \bkR{ 
			\gfcTZV{l}{x_{\sig(2)\sig(3)},x_{\sig(2)\sig(3)},x_{\sig(3)}} 
			+ 
			\gfcTZV{l}{x_{\sig(2)\sig(3)},x_{\sig(3)},x_{\sig(3)}} 
			+ 
			\gfcTZV{l}{x_{\sig(2)\sig(3)},x_{\sig(3)},0} 
		}
		\nonumber
	}
	{
		-
		\Sm{\sig\in\lsetS{32} } \bkR{ 
			\gfcTZV{l}{x_{\sig(2)},x_{\sig(2)},x_{\sig(3)}} 
			+ 
			\gfcTZV{l}{x_{\sig(2)},x_{\sig(3)},x_{\sig(3)}} 
			+ 
			\gfcTZV{l}{x_{\sig(2)},x_{\sig(3)},0} 
		}
		\lnAHs{10}
		+
		\opF{\gfcDZV{l}{x_{3},x_{2}}}{x_{3}} 
		+
		\opF{x_2^{l-2}-x_3^{l-2}}{x_2-x_3} \lmfcZ{1,l-1}
		-
		x_2^{l-3}\fcZ{l-2}\fcZ{2}
	}	
Subtracting \refEq{3_Thm1P_Eq2} from \refEq{3_Thm1P_Eq1} gives \refEq{1_Thm1i_Eq1}.

By \refEq{3_Prop1ii_Eq1}, we also have
	\envLinePd[3_Thm1P_Eq3]
	{
		\Sm{\sig\in\gpSym{3}} \bkS[G]{
			\gfcTZV{l}{x_{\sig(1)\sig(2)\sig(3)},x_{\sig(2)\sig(3)},x_{\sig(3)}}
			\lnAHs[]{10}
			-
			\bkR[B]{
				\gfcTZV{l}{x_{\sig(2)\sig(3)},x_{\sig(2)\sig(3)},x_{\sig(3)}}
				+
				\gfcTZV{l}{x_{\sig(2)\sig(3)},x_{\sig(3)},x_{\sig(3)}}
				+
				\gfcTZV{l}{x_{\sig(2)\sig(3)},x_{\sig(3)},0}
			}
		}
		\nonumber
	}
	{
		\Sm{\sig\in\gpSym{3}}
		\bkS[G]{
			\gfcTZV{l}{x_{\sig(1)},x_{\sig(2)},x_{\sig(3)}}
			\lnAHs{25}
			-
			\bkR[B]{ \gfcTZV{l}{x_{\sig(2)},x_{\sig(2)},x_{\sig(3)}} + \gfcTZV{l}{x_{\sig(2)},x_{\sig(3)},x_{\sig(3)}} + \gfcTZV{l}{x_{\sig(2)},x_{\sig(3)},0} }
			\lnAHs{200}
			+
			\opF{\gfcDZV{l}{x_{\sig(1)},x_{\sig(2)}}}{x_{\sig(1)} - \pw{x_{\sig(3)}}{1}} 
			+
			\opF{\gfcDZV{l}{x_{\sig(1)},x_{\sig(2)}}}{x_{\sig(2)} - \pw{x_{\sig(3)}}{1}} 
		}
		\lnAHs{0}
		+
		\Sm{\tau\in\gpAlt{3}}
		\bkS[G]{
			\opF{x_{\tau(2)}^{l-2}-x_{\tau(3)}^{l-2}}{x_{\tau(2)}-x_{\tau(3)}} \lmfcZ{1,l-1}
			+
			\bkR[g]{x_{\tau(1)}^{l-1} \PdT{i=2}{3} \opF{1}{x_{\tau(1)}-x_{\tau(i)}}} \fcZ{l}
			-
			x_{\tau(1)}^{l-3} \fcZ{l-2}\fcZ{2}
		}
	}
Since 
	\envM{\gpAlt{3}}{ \bkB{\lue, (123), (132)} } and 
	\envMCm
	{
		\gpSym{3}
	}
	{
		\bkB{ \tau\sig \mVert \tau\in\gpAlt{3}, \sig\in\lsetS{32} }
	}
	taking the summation of \refEq{3_Thm1P_Eq2} with $(x_2,x_3)=(x_{\tau(2)},x_{\tau(3)})$ over all $\tau\in\gpAlt{3}$ yields
	\envLinePd
	{
		\Sm{\sig\in\gpSym{3} } \bkR{ 
			\gfcTZV{l}{x_{\sig(2)\sig(3)},x_{\sig(2)\sig(3)},x_{\sig(3)}} 
			+ 
			\gfcTZV{l}{x_{\sig(2)\sig(3)},x_{\sig(3)},x_{\sig(3)}} 
			+ 
			\gfcTZV{l}{x_{\sig(2)\sig(3)},x_{\sig(3)},0} 
		}
		\lnAHs{0}
		+
		\Sm{\tau\in\gpAlt{3}} \opF{\gfcDZV{l}{x_{\tau(3)},x_{\tau(2)}}}{x_{\tau(3)}} 
	}
	{
		\Sm{\sig\in\gpSym{3} } \bkR{ 
			\gfcTZV{l}{x_{\sig(2)},x_{\sig(2)},x_{\sig(3)}} 
			+ 
			\gfcTZV{l}{x_{\sig(2)},x_{\sig(3)},x_{\sig(3)}} 
			+ 
			\gfcTZV{l}{x_{\sig(2)},x_{\sig(3)},0} 
		}
		\lnAHs{0}
		+
		\Sm{\tau\in\gpAlt{3}} \bkS[G]{ 
			\gfcTZV{l}{x_{\tau(3)},x_{\tau(2)\tau(3)},x_{\tau(3)}}+\gfcTZV{l}{x_{\tau(3)},x_{\tau(3)\tau(2)},x_{\tau(2)}} 
			\lnAHs{150}
			-
			\opF{x_{\tau(2)}^{l-2}-x_{\tau(3)}^{l-1}}{x_{\tau(2)}-x_{\tau(3)}} \lmfcZ{1,l-1}
			+
			x_{\tau(1)}^{l-3}\fcZ{l-2}\fcZ{2}
		}
	}
By adding this to \refEq{3_Thm1P_Eq3}, we obtain \refEq{1_Thm1ii_Eq1}.
}
%Z
%A
%:3_Rem1
\begin{remark}\label{3_Rem1}
As mentioned at the beginning of this section,
	the coefficients of the $x_1^{p-1}x_2^{q-1}x_3^{r-1}$ terms of \refEq{3_Prop1i_Eq1} 
	are equal to the extended double shuffle relations $\lclsEDS{l}{2,1}$.
Thus those of \refEq{3_Thm1P_Eq1} are also equal to $\lclsEDS{l}{2,1}$
	because \refEq{3_Thm1P_Eq1} is a restatement of \refEq{3_Prop1i_Eq1},
	but 
	it is difficult to claim that \refEq{3_Thm1P_Eq1} is beautiful. 
In contrast, 
	those of \refEq{1_Thm1i_Eq1} are not equal to $\lclsEDS{l}{2,1}$,
	because \refEq{1_Thm1i_Eq1} is obtained by removing the coefficients of the $x_2^{q-1}x_3^{r-1}$ terms of \refEq{3_Thm1P_Eq1},
	but \refEq{1_Thm1i_Eq1} is more beautiful than \refEq{3_Thm1P_Eq1}.
\end{remark}%Z

\section{\sectFour} \label{sectFour}
In this section,	
	we prove \refThm{1_Thm2}, which is the parameterized sum formula \refEq{1_Thm2_Eq1} for triple zeta values,
	and also give some weighted sum formulas as a corollary of the theorem.
%A
\envProof[\refThm{1_Thm2}]{
We easily see that
	\envHLineCFPd
	{
		\gpAlt{3}
	}
	{
		\bkB[n]{ \tau\sig \mVert \tau\in\gpAlt{3}, \sig=(321) }
	}
	{
		\gpSym{3}
	}
	{
		\bkB[n]{ \tau\sig \mVert \tau\in\gpAlt{3}, \sig\in\lsetS{31} }
	\lnP{=}
		\bkB[n]{ \tau\sig \mVert \tau\in\gpAlt{3}, \sig\in\lsetS{32} }
	}
Since 
	\envM{ \lsetS{321} }{ \bkB{(321)} \cup \lsetS{32} },
	taking the summation of \refEq{1_Thm1i_Eq1} with $(x_1,x_2,x_3)=(x_{\tau(1)},x_{\tau(2)},x_{\tau(3)})$ over all $\tau\in\gpAlt{3}$ gives
	\envLinePd
	{
		\Sm{\tau\in\gpAlt{3}} 
		\bkS[g]{
			\gfcTZV{l}{ x_{\tau(1)\tau(3)},x_{\tau(2)\tau(3)},x_{\tau(3)}}
			+
			\gfcTZV{l}{ x_{\tau(1)\tau(3)},x_{\tau(3)\tau(2)},x_{\tau(2)}}
			\lnAHs{120}
			+
			\gfcTZV{l}{x_{\tau(3)\tau(1)},x_{\tau(1)},x_{\tau(2)}}
			+
			\opF{\gfcDZV{l}{x_{\tau(3)},x_{\tau(2)}}}{x_{\tau(3)}} 
		}
	}
	{
		\Sm{\sig\in\gpSym{3}} \gfcTZV{l}{x_{\sig(1)},x_{\sig(2)},x_{\sig(3)}} 
		\lnAHs{0}
		+
		\Sm{\tau\in\gpAlt{3}} \bkR{ 
			\gfcTZV{l}{x_{\tau(1)},x_{\tau(2)},x_{\tau(3)}}
			+
			\gfcTZV{l}{x_{\tau(3)},x_{\tau(2)\tau(3)},x_{\tau(3)}} 
			+ 
			\gfcTZV{l}{x_{\tau(3)},x_{\tau(3)\tau(2)},x_{\tau(2)}} 
		}
		\lnAHs{0}
		+
		\Sm{\sig\in\gpSym{3}}  \bkR{ 
			\opF{\gfcDZV{l}{x_{\sig(1)},x_{\sig(2)}}}{x_{\sig(1)} - \pw{x_{\sig(3)}}{1}} 
			+
			\opF{\gfcDZV{l}{x_{\sig(1)},x_{\sig(2)}}}{x_{\sig(2)} - \pw{x_{\sig(3)}}{1}} 	 
		}
	}
Subtracting this from \refEq{1_Thm1ii_Eq1} yields \refEq{1_Thm2_Eq1},
	and \refEq{1_Thm2_Eq1} with \refEq{1_Plane_DefGfcTZV} and  \refEq{2.1_Prop1PEq2} proves \refEq{1_Thm2_Eq2}.
}
%Z

We give some weighted sum formulas in \refEq{4_Cor1_Eq1} below for triple zeta values by substituting appropriate values for 
	$(x_1,x_2,x_3)$ in formula \refEq{1_Thm2_Eq2}.
The second equation of \refEq{4_Cor1_Eq1} is
	the result of Guo and Xie \cite[Theorem 1.1]{GX09} in the case of triple zeta values.
%A
%:4_Cor1
\begin{corollary}[\text{cf. \cite[Theorem 1.1]{GX09}}]\label{4_Cor1}
Let $l$ be a positive integer such that $l\geq4$. 
Then we have
	\envHLineCSPd[4_Cor1_Eq1]
	{
		\tpSm{l_1\geq2,l_2,l_3\geq1}{l_1+l_2+l_3=l} 
		\bkR{ \pw{3}{l_1-1}\pw{2}{l_2}-\pw{2}{l_1+l_2-1}-\pw{2}{l_1-1} }\fcZ{l_1,l_2,l_3}
	}
	{
		\opF{(l-4)(l+1)}{6}\fcZ{l}
	}
	{
		\tpSm{l_1\geq2,l_2,l_3\geq1}{l_1+l_2+l_3=l} 
		\bkR{ \pw{2}{l_1+l_2-1}+\pw{2}{l_1-1}-\pw{2}{l_2} }\fcZ{l_1,l_2,l_3}
	}
	{
		l\fcZ{l}
	}
	{
		\tpSm{l_1\geq2,l_2,l_3\geq1}{l_1+l_2+l_3=l} 
		\bkR{ \pw{3}{l_1-1} - 1 } \pw{2}{l_2} \fcZ{l_1,l_2,l_3}
	}
	{
		\opF{(l-1)(l+4)}{6}\fcZ{l}
	}
\end{corollary}
%Z
%A
\envProof{
By substituting $(1,1,1)$ for $(x_1,x_2,x_3)$ in \refEq{1_Thm2_Eq2} and using \refEq{1_Plane_EqSFofTZV},
	we obtain the first equation of \refEq{4_Cor1_Eq1}
	since $\tpSm{l_1,l_2,l_3\geq1}{l_1+l_2+l_3=l} 1 = (l-2)(l-1)/2$.
By substituting $(1,1,0)$ for $(x_1,x_2,x_3)$ in \refEq{1_Thm2_Eq2} and using \refEq{1_Plane_EqSFofTZV},
	we obtain the second equation.
The third equation is obtained by adding the first to the second.
}
%Z

\section{\sectFive} \label{sectFive}
In this final section,
	we prove \refThm{1_Thm3}, which gives the restricted sum formulas \refEq{1_Thm3i_Eq1} and \refEq{1_Thm3ii_Eq1}.
We also present some restricted analogues of the formulas given by Granville, Hoffman, and Ohno, as \refProp{5_Prop1}.

We prepare \refLem[s]{5_Lem1} and \ref{5_Lem2} to prove the theorem and proposition.
\refLem{5_Lem1} states that
	the restricted sums of triple zeta values of the forms
	\envPLine
	{
		\tpThSm{l_1\geq2,l_2,l_3\geq1}{l_1+l_2+l_3=l}{l_i:\emph{even or odd}} \fcZ{l_1,l_2,l_3},
		\qquad
		\tpThSm{l_1\geq2,l_3\geq1}{l_1+l_3=l-1}{l_i:\emph{even or odd}} \fcZ{l_1,1,l_3},
		\qquad
		\tpThSm{l_1\geq2,l_2\geq1}{l_1+l_2=l-1}{l_i:\emph{even or odd}} \fcZ{l_1,l_2,1}
	}
	can be written in terms of combinations of $\gfcTZV{l}{\ep_1,\ep_2,\ep_3}$'s where $\ep_1,\ep_2,\ep_3\in\bkB{0,\pm1}$,
	and \refLem{5_Lem2} evaluates certain combinations of $\gfcTZV{l}{\ep_1,\ep_2,\ep_3}$'s.
Note that
	the first equation of \refEq{5_Lem2i_Eq1} in \refLem{5_Lem2} has already been proved by Shen and Cai (see \cite[(4)]{SC12}).
We often use the following equations, which are obvious by the definitions, in the proofs without comment:
	\envHLineCFPd[5_Plane_EqOfDZVandTZV]
	{
		\gfcDZV{l}{-x_1,-x_2}
	}
	{
		\mo^{l}\gfcDZV{l}{x_1,x_2}
	}
	{
		\gfcTZV{l}{-x_1,-x_2,-x_3}
	}
	{
		\mo^{l-1}\gfcTZV{l}{x_1,x_2,x_3}
	}

%A
%:5_Lem1
\begin{lemma}\label{5_Lem1}
Let $l$ be an integer such that $l\geq4$, 
	and $\pSm[t]{P(l_1,\ldots,l_n)}$ be as in \refThm{1_Thm3}.
\\{\bf (i)} 
For $\alp_1,\alp_2,\alp_3,\alp_4\in\bkB{\pm1}$, 
	we define 
	\envHLineDefPd
	{
		\lfcST{\alp_1,\alp_2,\alp_3,\alp_4}
	}
	{
		\alp_1\gfcTZV{l}{1,1,1} + \alp_2\gfcTZV{l}{1,1,-1} + \alp_3\gfcTZV{l}{1,-1,1} + \alp_4\gfcTZV{l}{1,-1,-1}
	}
Then we have
	\envHLineCE[5_Lem1_Eq1]
	{
		\opF{ \lfcST{1,-1,-1,1} }{4}
	}
	{
		\envCaseT[d]{	
			\pSm{l_1,l_2,l_3:even} \fcZ{l_1,l_2,l_3} 				&	(l:even)
		}{
			\pSm{l_2,l_3:even \atop l_1:odd} \fcZ{l_1,l_2,l_3} 		&	(l:odd),
		}
	}
	{
		\opF{ \lfcST{1,-1,1,-1} }{4}
	}
	{
		\envCaseT[d]{	
			\pSm{l_3:even \atop l_1,l_2:odd} \fcZ{l_1,l_2,l_3}		&	(l:even)
		}{
			\pSm{l_1,l_3:even \atop l_2:odd} \fcZ{l_1,l_2,l_3}		&	(l:odd),
		}
	}
	{
		\opF{ \lfcST{1,1,-1,-1} }{4}
	}
	{
		\envCaseT[d]{	
			\pSm{l_2:even \atop l_1,l_3:odd} \fcZ{l_1,l_2,l_3}		&	(l:even)
		}{
			\pSm{l_1,l_2:even \atop l_3:odd} \fcZ{l_1,l_2,l_3}		&	(l:odd),
		}
	}
	{
		\opF{ \lfcST{1,1,1,1} }{4}
	}
	{
		\envCaseT[d]{	
			\pSm{l_1:even \atop l_2,l_3:odd} \fcZ{l_1,l_2,l_3}		&	(l:even)
		}{
			\pSm{l_1,l_2,l_3:odd}	 \fcZ{l_1,l_2,l_3}	&	(l:odd).
		}
	}
{\bf (ii)} We have
	\envHLineCE[5_Lem1_Eq2]
	{
		\opF{ \gfcTZV{l}{1,0,1}-\gfcTZV{l}{-1,0,1} }{2}
	}
	{
		\envCaseT[d]{	
			\pSm{l_1:even, l_3:odd \atop l_2=1} \fcZ{l_1,1,l_3}			&	(l:even)
		}{
			\pSm{l_1,l_3: even \atop l_2=1} \fcZ{l_1,1,l_3} 				&	(l:odd),
		}
	}
	{
		\opF{ \gfcTZV{l}{1,0,1}+\gfcTZV{l}{-1,0,1} }{2}
	}
	{
		\envCaseT[d]{	
			\pSm{l_1:odd, l_3:even \atop l_2=1} \fcZ{l_1,1,l_3}			&	(l:even)
		}{
			\pSm{l_1,l_3: odd \atop l_2=1} \fcZ{l_1,1,l_3}				&	(l:odd),
		}
	}
	{
		\opF{ \gfcTZV{l}{1,1,0}-\gfcTZV{l}{-1,1,0} }{2}
	}
	{
		\envCaseT[d]{	
			\pSm{l_1:even, l_2:odd \atop l_3=1} \fcZ{l_1,l_2,1}			&	(l:even)
		}{
			\pSm{l_1,l_2: even \atop l_3=1} \fcZ{l_1,l_2,1} 				&	(l:odd),
		}
	}
	{
		\opF{ \gfcTZV{l}{1,1,0}+\gfcTZV{l}{-1,1,0} }{2}
	}
	{
		\envCaseT[d]{	
			\pSm{l_1:odd, l_2:even \atop l_3=1} \fcZ{l_1,l_2,1}			&	(l:even)
		}{
			\pSm{l_1,l_2: odd \atop l_3=1} \fcZ{l_1,l_2,1}				&	(l:odd).
		}
	}
\end{lemma}
%Z
%A
\envProof{
For $\ep_1,\ep_2\in\bkB{\pm1}$, we have
	\envHLineCFPd
	{
		\lfcSTid{-}{\ep_1,\ep_2}
	\lnP{:=}
		\gfcTZV{l}{\ep_1,\ep_2,1} - \gfcTZV{l}{\ep_1,\ep_2,-1}
	}
	{
		2 \pSm{l_3:even} \ep_1^{l_1-1}\ep_2^{l_2-1} \fcZ{l_1,l_2,l_3}
	}
	{
		\lfcSTid{+}{\ep_1,\ep_2}
	\lnP{:=}
		\gfcTZV{l}{\ep_1,\ep_2,1} + \gfcTZV{l}{\ep_1,\ep_2,-1}
	}
	{
		2 \pSm{l_3:odd} \ep_1^{l_1-1}\ep_2^{l_2-1} \fcZ{l_1,l_2,l_3}
	}
Then, for $\ep_1\in\bkB{\pm1}$, we obtain
	\envHLineCEPd
	{
		\lfcSTid{-}{\ep_1,1} - \lfcSTid{-}{\ep_1,-1}
	}
	{
		4 \pSm{l_2,l_3:even} \ep_1^{l_1-1} \fcZ{l_1,l_2,l_3}
	}
	{
		\lfcSTid{-}{\ep_1,1} + \lfcSTid{-}{\ep_1,-1}
	}
	{
		4 \pSm{l_3:even \atop l_2:odd} \ep_1^{l_1-1} \fcZ{l_1,l_2,l_3}
	}
	{
		\lfcSTid{+}{\ep_1,1} - \lfcSTid{+}{\ep_1,-1}
	}
	{
		4 \pSm{l_2:even \atop l_3:odd } \ep_1^{l_1-1} \fcZ{l_1,l_2,l_3}
	}
	{
		\lfcSTid{+}{\ep_1,1} + \lfcSTid{+}{\ep_1,-1}
	}
	{
		4 \pSm{l_2,l_3:odd } \ep_1^{l_1-1} \fcZ{l_1,l_2,l_3}
	}
Thus we conclude that
	\envPLine{\HLineCECm[p]
	{
		\bkR[b]{ \lfcSTid{-}{1,1} - \lfcSTid{-}{1,-1} } - \bkR[b]{ \lfcSTid{-}{-1,1} - \lfcSTid{-}{-1,-1} }
	}
	{
		8 \pSm{l_1,l_2,l_3:even} \fcZ{l_1,l_2,l_3}
	}
	{
		\bkR[b]{ \lfcSTid{-}{1,1} - \lfcSTid{-}{1,-1} } + \bkR[b]{ \lfcSTid{-}{-1,1} - \lfcSTid{-}{-1,-1} }
	}
	{
		8 \pSm{l_2,l_3:even \atop l_1:odd } \fcZ{l_1,l_2,l_3}
	}
	{
		\bkR[b]{ \lfcSTid{-}{1,1} + \lfcSTid{-}{1,-1} } - \bkR[b]{ \lfcSTid{-}{-1,1} + \lfcSTid{-}{-1,-1} }
	}
	{
		8 \pSm{l_1,l_3:even \atop l_2:odd } \fcZ{l_1,l_2,l_3}
	}
	{
		\bkR[b]{ \lfcSTid{-}{1,1} + \lfcSTid{-}{1,-1} } + \bkR[b]{ \lfcSTid{-}{-1,1} + \lfcSTid{-}{-1,-1} }
	}
	{
		8 \pSm{l_3:even \atop l_1,l_2:odd } \fcZ{l_1,l_2,l_3}
	}\HLineCEPd
	{
		\bkR[b]{ \lfcSTid{+}{1,1} - \lfcSTid{+}{1,-1} } - \bkR[b]{ \lfcSTid{+}{-1,1} - \lfcSTid{+}{-1,-1} }
	}
	{
		8 \pSm{l_1,l_2:even \atop l_3:odd } \fcZ{l_1,l_2,l_3}
	}
	{
		\bkR[b]{ \lfcSTid{+}{1,1} - \lfcSTid{+}{1,-1} } + \bkR[b]{ \lfcSTid{+}{-1,1} - \lfcSTid{+}{-1,-1} }
	}
	{
		8 \pSm{l_2:even \atop l_1,l_3:odd } \fcZ{l_1,l_2,l_3}
	}
	{
		\bkR[b]{ \lfcSTid{+}{1,1} + \lfcSTid{+}{1,-1} } - \bkR[b]{ \lfcSTid{+}{-1,1} + \lfcSTid{+}{-1,-1} }
	}
	{
		8 \pSm{l_1:even \atop l_2,l_3:odd } \fcZ{l_1,l_2,l_3}
	}
	{
		\bkR[b]{ \lfcSTid{+}{1,1} + \lfcSTid{+}{1,-1} } + \bkR[b]{ \lfcSTid{+}{-1,1} + \lfcSTid{+}{-1,-1} }
	}
	{
		8 \pSm{l_1,l_2,l_3:odd } \fcZ{l_1,l_2,l_3}
	}
	}
These equations yield \refEq{5_Lem1_Eq1}
	since we see from \refEq{5_Plane_EqOfDZVandTZV} that 
	\envHLineCFPd
	{
		\lfcSTid{-}{-\ep_1,-\ep_2}
	}
	{
		\mo^{l} \lfcSTid{-}{\ep_1,\ep_2}
	}
	{
		\lfcSTid{+}{-\ep_1,-\ep_2}
	}
	{
		\mo^{l-1} \lfcSTid{+}{\ep_1,\ep_2}
	}
	
We can prove \refEq{5_Lem1_Eq2} similarly and more easily, and so omit the proof here.
}
%Z
%A
%:5_Lem2
\begin{lemma}[\text{cf. \cite[(4)]{SC12}}]\label{5_Lem2}
Let $l$ be an integer such that $l\geq4$. 
\mbox{}\\{\bf (i)} If $l$ is even, then
	\envHLineCSPd[5_Lem2i_Eq1]
	{
		 \gfcTZV{l}{1,1,-1} + \gfcTZV{l}{1,-1,1} + \gfcTZV{l}{-1,1,1}
	}
	{
		-\opF{3}{2}\fcZ{l} + \fcZ{l-2}\fcZ{2}
	}
	{
		\gfcTZV{l}{1,1,-1} - \gfcTZV{l}{-1,1,1}
	}
	{
		-\opF{1}{2}\fcZ{l} + \fcZ{l-2}\fcZ{2}
	}
	{
		\gfcTZV{l}{-1,0,1}
	}
	{
		-\opF{1}{2}\fcZ{l} +\fcZ{l-1,1}-\fcZ{l-2,2}
	}
{\bf (ii)} If $l$ is odd, then
	\envHLineCEPd[5_Lem2ii_Eq1]
	{
		\gfcTZV{l}{1,1,-1} + \gfcTZV{l}{1,-1,1} + \gfcTZV{l}{-1,1,1} + \gfcDZV{l}{-1,1}
	}
	{
		- \opF{1}{2}\fcZ{l} - \fcZ{2,l-2} 
	}
	{
		\gfcTZV{l}{1,1,-1} + \gfcTZV{l}{-1,1,1}
	}
	{
		- \opF{1}{2}\fcZ{l} - \fcZ{2,l-2} 
	}
	{
		\gfcTZV{l}{-1,0,1}
	}
	{
		-\opF{1}{2}\fcZ{l} - \fcZ{l-1,1}
	}
	{
		\gfcTZV{l}{-1,1,0} + \opF{1}{2}\gfcDZV{l}{-1,1}
	}
	{
		-\opF{1}{4}\fcZ{l}
		\lnAHs{0} 
		- 
		\opF{1}{2}\fcZ{l-1,1} - \opF{1}{2}\fcZ{l-2,2}
	}
\end{lemma}
%Z
%A
\envProof{
We obtain the following equations \refEq{5_Lem2P_Eq1a} and \refEq{5_Lem2P_Eq1b} 
	from \refEq{1_Thm1i_Eq1} with $(x_1,x_2,x_3)=(-1,0,1)$ and $(-1,1,0)$, respectively,
	and \refEq{5_Lem2P_Eq1c} from \refEq{1_Thm2_Eq1} with $(-1,0,1)$;
	\envHLineCSCmNme
	{\label{5_Lem2P_Eq1a}
		\gfcTZV{l}{-1,0,1} + \bkR[b]{1-\mo^l}\gfcTZV{l}{-1,1,0} + \gfcDZV{l}{-1,1}
	}
	{
		-
		\fcZ{l} 
		\lnAHs[]{0}
		- 
		\opF{1-3\mo^l}{2}\fcZ{l-1,1} 
		\lnAHs{0}
		- 
		\fcZ{l-2,2}
		\nonumber
	}
	{\label{5_Lem2P_Eq1b}
		\gfcTZV{l}{-1,1,1} + \gfcTZV{l}{-1,-1,1} - \gfcTZV{l}{-1,0,1}
	}
	{
		-\mo^l \fcZ{l-1,1}
		\lnAHs[]{0} 
		- 
		\fcZ{2,l-2}
		\nonumber
	}
	{\label{5_Lem2P_Eq1c}
		\gfcTZV{l}{-1,1,1}+ \gfcTZV{l}{-1,-1,1} + \gfcTZV{l}{1,0,-1}
	}
	{
		-\opF{1-\mo^l}{2}\fcZ{l} 
		\lnAHs[]{5}
		- 
		\fcZ{l-1,1} - \fcZ{2,l-2}
		\nonumber
	}
	where we use \refEq{1_Plane_EqGranville2} and \refEq{1_Plane_EqGranville1} below 
	for the proofs of \refEq{5_Lem2P_Eq1a} and \refEq{5_Lem2P_Eq1c}, respectively.
By taking the difference of \refEq{5_Lem2P_Eq1b} and \refEq{5_Lem2P_Eq1c}, we also obtain
	\envHLinePd[5_Lem2P_Eq2]
	{
		\bkR{ 1- \mo^{l} } \gfcTZV{l}{-1,0,1}
	}
	{
		-\opF{1-\mo^l}{2}\fcZ{l} - \bkR{ 1 - \mo^l }\fcZ{l-1,1} 
	}
	
The third equation of \refEq{5_Lem2i_Eq1} follows from \refEq{5_Lem2P_Eq1a} and \refEq{1_Plane_EqRSFofDZV},
	and the third equation of \refEq{5_Lem2ii_Eq1} follows from \refEq{5_Lem2P_Eq2}.
(Note the accompanying assumptions; that is, $l$ is even if we are considering \refEq{5_Lem2i_Eq1} and odd if we are considering \refEq{5_Lem2ii_Eq1}.)
The fourth equation of \refEq{5_Lem2ii_Eq1} is derived from the third equation and \refEq{5_Lem2P_Eq1a}.

The second equation of \refEq{5_Lem2i_Eq1} 
	follows from the third equation and \refEq{5_Lem2P_Eq1b},
	and
	the first equation of \refEq{5_Lem2i_Eq1} is also derived from the third equation and 
	\envLineCm
	{
		\gfcTZV{l}{1,1,-1}+\gfcTZV{l}{1,-1,1}+\gfcTZV{l}{-1,1,1}+\gfcTZV{l}{-1,0,1}+\opF{1-\mo^l}{2}\gfcDZV{l}{-1,1}
	}
	{
		-\fcZ{l} +\mo^{l}\fcZ{l-1,1}+\mo^{l}\fcZ{2,l-2}
	}
	which follows from \refEq{1_Thm1i_Eq1} with $(x_1,x_2,x_3)=(-1,-1,1)$ and \refEq{1_Plane_EqGranville1}.
The first and second equations of \refEq{5_Lem2ii_Eq1} follow similarly as those of \refEq{5_Lem2i_Eq1}.
}%Z
We prove \refThm{1_Thm3}.
%A
\envProof{
First, we consider the case that $l$ is even.
We see from \refEq{5_Plane_EqOfDZVandTZV}  and \refEq{5_Lem1_Eq1} that
	\envHLineCFPd
	{
		\pSm{l_1,l_2,l_3:even} \fcZ{l_1,l_2,l_3} 
	}
	{
		\opF{\gfcTZV{l}{1,1,1}-\gfcTZV{l}{1,1,-1}-\gfcTZV{l}{1,-1,1}-\gfcTZV{l}{-1,1,1}}{4}
	}
	{
		\bkR{ \rTx{ \pSm{l_1:even \atop l_2,l_3:odd} - \pSm{l_3:even \atop l_1,l_2:odd} } } \fcZ{l_1,l_2,l_3}
	}
	{
		\opF{\gfcTZV{l}{1,1,-1}-\gfcTZV{l}{-1,1,1}}{2}
	}
These together with \refEq{1_Plane_EqSFofTZV} and \refEq{5_Lem2i_Eq1} yield the first and third equations of \refEq{1_Thm3i_Eq1}.
The second equation of \refEq{1_Thm3i_Eq1} is easily derived from the first equation and \refEq{1_Plane_EqSFofTZV}.

Next we consider the case that $l$ is odd.
We similarly see that
	\envHLineCFPd[5_Thm3P_Eq1]
	{
		\pSm{l_1,l_2,l_3:odd} \fcZ{l_1,l_2,l_3}
	}
	{
		\opF{\gfcTZV{l}{1,1,1}+\gfcTZV{l}{1,1,-1}+\gfcTZV{l}{1,-1,1}+\gfcTZV{l}{-1,1,1}}{4}
	}
	{
		\pSm{l_1,l_3:even \atop l_2:odd} \fcZ{l_1,l_2,l_3}
	}
	{
		\opF{\gfcTZV{l}{1,1,1}-\gfcTZV{l}{1,1,-1}+\gfcTZV{l}{1,-1,1}-\gfcTZV{l}{-1,1,1}}{4}
	}
The first equation of \refEq{1_Thm3ii_Eq1} follows from the first equations of \refEq{5_Lem2ii_Eq1} and \refEq{5_Thm3P_Eq1}.
From the difference of the first and second equations of \refEq{5_Lem2ii_Eq1}, we obtain
	\envMPd
	{
		\gfcTZV{l}{1,-1,1}
	}
	{
		- \gfcDZV{l}{-1,1}
	}
Thus, the second equation of \refEq{1_Thm3ii_Eq1} follows from the second equations of \refEq{5_Lem2ii_Eq1} and \refEq{5_Thm3P_Eq1}.
The third equation of \refEq{1_Thm3ii_Eq1} is easily derived from 
	the first and second together with \refEq{1_Plane_EqSFofDZV} and \refEq{1_Plane_EqSFofTZV}.
}%Z

Finally, we give 
	restricted analogues of the following formulas given by Granville, Hoffman and Ohno: 
	\envHLineCFNmePd
	{\label{1_Plane_EqGranville1}
		\gfcTZV{l}{1,0,1}
	\lnP{=}
		\tpSm{l_1\geq2,l_3\geq1}{l_1+l_3=l-1} \fcZ{l_1,1,l_3}
	}
	{
		\fcZ{l-1,1} + \fcZ{2,l-2}
	}
	{\label{1_Plane_EqGranville2}
		\gfcTZV{l}{1,1,0}
	\lnP{=}
		\tpSm{l_1\geq2,l_2\geq1}{l_1+l_2=l-1} \fcZ{l_1,l_2,1}
	}
	{
		\fcZ{l-1,1} + \fcZ{l-2,2}
	}
(Refer to \cite[(9)]{Granville97} and \cite[(1)]{HO03} with $(k_1,k_2)=(l-2,1)$ for \refEq{1_Plane_EqGranville1}, 
	and \cite[(10)]{Granville97} and \cite[Theorem 5.1]{Hoffman92} with $(i_1,i_2)=(l-2,1)$ for \refEq{1_Plane_EqGranville2}.)
We note that
	\refEq{1_Plane_EqGranville1} and \refEq{1_Plane_EqGranville2}
	are directly derived from \refEq{1_Thm1i_Eq1} with $(x_1, x_2, x_3)=(1,1,0)$ and $(1,0,0)$, respectively.
%A
%:5_Prop1
\begin{proposition}\label{5_Prop1}
Let $l$ be an integer such that $l\geq4$, and $\pSm[t]{P(l_1,\ldots,l_n)}$ be as in \refThm{1_Thm3}.
\mbox{}\\{\bf (i)} If $l$ is even, then
\envHLineCFPd[5_Prop1i_Eq1]
	{
		\pSm{l_1:even, l_3:odd \atop l_2=1} \fcZ{l_1,1,l_3}
	}
	{
		\opF{1}{4}\fcZ{l} + \opF{1}{2}\fcZ{l-2,2} + \opF{1}{2}\fcZ{2,l-2}
	}
	{
		\pSm{l_1:odd, l_3:even \atop l_2=1} \fcZ{l_1,1,l_3}
	}
	{
		-\opF{1}{4}\fcZ{l} + \fcZ{l-1,1} - \opF{1}{2} \fcZ{l-2,2} + \opF{1}{2}\fcZ{2,l-2}
	}
{\bf (ii)} If $l$ is odd, then
\envHLineCEPd[5_Prop1ii_Eq1]
	{
		\pSm{l_1,l_3: even \atop l_2=1} \fcZ{l_1,1,l_3}
	}
	{
		\opF{1}{4}\fcZ{l} + \fcZ{l-1,1} + \opF{1}{2} \fcZ{2,l-2}
	}
	{
		\pSm{l_1,l_3: odd \atop l_2=1} \fcZ{l_1,1,l_3}
	}
	{
		-\opF{1}{4}\fcZ{l} + \opF{1}{2} \fcZ{2,l-2}
	}
	{
		\pSm{l_1,l_2: even \atop l_3=1} \fcZ{l_1,l_2,1} + \opF{1}{2}\pSm{l_1:even \atop l_2:odd}\fcZ{l_1,l_2}
	}
	{
		\opF{3}{8}\fcZ{l} + \opF{3}{4} \fcZ{l-1,1} + \opF{3}{4} \fcZ{l-2,2}
	}
	{
		\pSm{l_1,l_2: odd \atop l_3=1} \fcZ{l_1,l_2,1} + \opF{1}{2}\pSm{l_2:even \atop l_1:odd}\fcZ{l_1,l_2}
	}
	{
		\opF{1}{8}\fcZ{l} + \opF{1}{4} \fcZ{l-1,1} + \opF{1}{4} \fcZ{l-2,2}
	}
\end{proposition}
%Z
%A
\envProof{
By virtue of the first equation of \refEq{5_Lem1_Eq2} and \refEq{1_Plane_EqGranville1},
	the first equations of \refEq{5_Prop1i_Eq1} and \refEq{5_Prop1ii_Eq1} follow from the third equations of \refEq{5_Lem2i_Eq1} and \refEq{5_Lem2ii_Eq1}, respectively.	
The second equations of \refEq{5_Prop1i_Eq1} and \refEq{5_Prop1ii_Eq1} are easily derived from the corresponding first equations together with \refEq{1_Plane_EqGranville1}.

Similarly,
	by virtue of the third equation of \refEq{5_Lem1_Eq2} and \refEq{1_Plane_EqGranville2},
 	the third equation of \refEq{5_Prop1ii_Eq1} follows from the fourth equation of \refEq{5_Lem2ii_Eq1}.
The fourth equation of \refEq{5_Prop1ii_Eq1} is easily derived from the third together with \refEq{1_Plane_EqSFofDZV} and \refEq{1_Plane_EqGranville2}.
}%Z

%--acknowledgement
\section*{Acknowledgements}
The author would like to thank Professors Kohji Matsumoto and Takashi Nakamura
	for introducing the papers \cite{KMT11} and \cite{SC12}, respectively.

%--reference

%--author information
\begin{flushleft}
\mbox{}\\ \qquad \mbox{}\\ \qquad
%\majorName		\mbox{}\\ \qquad
\departmentName	\mbox{}\\ \qquad
\organizationName	\mbox{}\\ \qquad
\placeAddress		\mbox{}\\ \qquad
\emailAddress
\end{flushleft}

\end{document}